\documentclass[a4paper]{jpconf}
\usepackage{graphicx}

\usepackage{amsfonts}
\usepackage{amssymb}
\usepackage{amsmath}
\usepackage{epic,color}
\usepackage{eepic}
\usepackage{epsf}

\def\blu#1{{\color{blue}#1}}
\def\red#1{{\color{red}#1}}

\def\be{\begin{equation}}
\def\eqn#1{\be\label{#1}}
\def\ee{\end{equation}}

\def\bea{\begin{eqnarray}}
\def\eqnn#1{\bea\label{#1}}
\def\eea{\end{eqnarray}}

\def\nn{\nonumber}
\def\nt{\noindent}
\def\nl{\hfill\break}

\def\nd{\end{document}}

\input epsf.tex
\newcount\figno
\def\fig#1#2#3{
\par\begingroup\parindent=0pt\leftskip=1cm\rightskip=1cm\parindent=0pt
\baselineskip=11pt \global\advance\figno by 1 
\epsfxsize=#3 \centerline{\epsfbox{#2}} \vskip 12pt #1\par
\endgroup\par}
\def\figlabel#1{\xdef#1{\the\figno}}
\def\encadremath#1{\vbox{\hrule\hbox{\vrule\kern8pt\vbox{\kern8pt
\hbox{$\displaystyle #1$}\kern8pt} \kern8pt\vrule}\hrule}}

  \def\tV{{\tilde V}}

\def\bu{\noindent $\bullet~$}
\def\rank{{\rm rank}}

\def\riga{-\kern-4pt - \kern-4pt -}
\font\fat=cmsy10 scaled\magstep5

\def\Bbullet{\raise-3pt\hbox{\fat\char"0F}}

\def\Box{
\vbox{ \halign to5pt{\strut##& \hfil ## \hfil \cr &$\kern -0.5pt
\sqcap$ \cr \noalign{\kern -5pt \hrule} }}~}

\def\down{\raise1.5pt\hbox{$\phantom{a}_2$}\downarrow}

\def\downa{\raise1.5pt\hbox{$\phantom{a}_{2\atop m_2}$}\downarrow}

\def\llr{\longrightarrow}

\def\({\left(}
\def\){\right)}

\def\lra{\longrightarrow}

\def\dia{{$\diamondsuit$}}

\def\ha{{\textstyle{\frac{1}{2}}}}  

\def\trha{{\textstyle{3\over2}}}

\def\bbc{\mathbb{C}}
\def\bac{\bbc} 
\def\bbr{\mathbb{R}}
\def\bbq{\mathbb{Q}}
\def\bbo{\mathbb{O}}
\def\bbn{\mathbb{N}}
\def\a{\alpha}
\def\b{\beta}

\def\vr{\vert}

\def\ed{\end{document}}

\def\ca{{\cal A}}  \def\cc{{\cal C}}
\def\cd{{\cal D}} \def\ce{{\cal E}} \def\cf{{\cal F}}
\def\cg{{\cal G}} \def\ch{{\cal H}} 
 \def\ck{{\cal K}} 
\def\cm{{\cal M}} \def\cn{{\cal N}} 
\def\cp{{\cal P}}  
 \def\ct{{\cal T}}

\def\ido{intertwining differential operator}
\def\idos{intertwining differential operators}

\def\L{\Lambda}
\def\r{\rho}

\begin{document}
 \pagestyle{empty}

\title{Multiplet classification for SU(n,n)}

\author{V.K.~Dobrev}

\address{Institute of Nuclear Research and Nuclear Energy,
 Bulgarian Academy of Sciences,
72 Tsarigradsko Chaussee, 1784 Sofia, Bulgaria}


\begin{abstract}
In the present paper we review our project of systematic
construction of invariant differential operators on the example of
the non-compact  algebras  $su(n,n)$ for $n=2,3,4$.
 We give explicitly the main multiplets of indecomposable elementary
representations   and some reduced multiplets.
We give explicitly the minimal representations.  Due to the recently
established parabolic relations
the multiplet classification results are valid also for the algebras $sl(2n,\bbr)$ and
when $n=2k$ for the algebras $su^*(4k)$ with
suitably chosen maximal parabolic subalgebras.
\end{abstract}

\section{Introduction}


 \blu{Invariant differential operators} play very important role in
the description of physical symmetries - starting from the early
occurrences in the Maxwell, d'Allembert, Dirac, equations,   to the
latest applications of (super-)differential operators in conformal
field theory, supergravity and string theory. Thus, it is important
for the applications in physics to study systematically such
operators.

In a recent paper \cite{Dobinv} we started the systematic explicit construction
of invariant differential operators. We gave an explicit description
of the building blocks, namely, the \blu{parabolic subgroups and
subalgebras} from which the necessary representations are induced.
Thus we have set the stage for study of different non-compact
groups.

Since the study and description of detailed classification should be
done group by group we had to decide which groups to study. Since
the most widely used algebras are the \blu{conformal algebras}
~\red{so(n,2)}~ in $n$-dimensional Minkowski space-time we
concentrated on a class that shares most of their properties.
  This class  consists of:
$$ so(n,2), ~~sp(n,\bbr), ~~su(n,n),  ~~so^*(4n),
~~E_{7(-25)}  $$   the corresponding analogs of Minkowski space-time
$V$ being:
$$\bbr^{n-1,1}, ~~{\rm Sym}(n,\bbr), ~~{\rm Herm}(n,\bbc),
~~{\rm Herm}(n,\bbq), ~~ {\rm Herm}(3,\bbo)$$
involving the four division algebras ~$\bbr,\bbc,\bbq,\bbo$.

 In view of applications to physics, we proposed to call these algebras   '\blu{conformal
Lie algebras}', (or groups)   \cite{Dobeseven}.

We have started the study of the above class in the framework of the
present approach in the cases: ~$so(n,2)$,  ~$su(n,n)$,
~$sp(n,\bbr)$, ~$E_{7(-25)}$,  cf. \cite{Dobeseven,Dobsunn,Dobspn,Dobparab}.

Lately, we discovered an efficient way to extend our considerations
beyond this class introducing the notion of 'parabolically related
non-compact semisimple Lie algebras' \cite{Dobparab}.

\bu {\it Definition:}  ~~~Let ~$\cg,\cg'$~ be two non-compact
semisimple Lie algebras with the same complexification ~$\cg^\bac
\cong \cg'^\bac$. We call them ~\red{parabolically related}~ if they
have parabolic subalgebras ~$\cp = \cm \oplus \ca \oplus \cn$,
~$\cp' = \cm' \oplus \ca' \oplus \cn'$, such that: ~$\cm^\bac
~\cong~ \cm'^\bac$~  ($\Rightarrow \cp^\bac ~\cong~ \cp'^\bac$).\dia

Certainly, there are many such parabolic relationships for any given
algebra ~$\cg$. Furthermore, two algebras ~$\cg,\cg'$~ may be
parabolically related with different parabolic subalgebras.

In the present paper we review our results on the case of
~$su(n,n)$, cf. \cite{Dobsunn,Dobsutt,Dobsuff}.
Due to the parabolic relationships these would be valid
also for ~$sl(2n,\bbr)$, and if ~$n=2k$~ also for ~$su^*(4k)$.

\section{Preliminaries}

 Let $G$ be a semisimple non-compact Lie group, and $K$ a
maximal compact subgroup of $G$. Then we have an {\it Iwasawa
decomposition} ~$G=KA_0N_0$, where ~$A_0$~ is Abelian simply
connected vector subgroup of ~$G$, ~$N_0$~ is a nilpotent simply
connected subgroup of ~$G$~ preserved by the action of ~$A_0$.
Further, let $M_0$ be the centralizer of $A_0$ in $K$. Then the
subgroup ~$P_0 ~=~ M_0 A_0 N_0$~ is a {\it minimal parabolic
subgroup} of $G$.  A {\it parabolic subgroup} ~$P ~=~ M' A' N'$~ is
any subgroup of $G$
which contains a minimal parabolic subgroup.

Further, let ~$\cg,\ck,\cp,\cm,\ca,\cn$~ denote the Lie algebras of
~$G,K,P,M,A,N$, resp.

For our purposes we need to restrict to  ~{\it maximal ~ parabolic
subgroups} ~$P=MAN$, i.e.  $\rank A =1$, resp. to ~{\it maximal ~
parabolic subalgebras} ~$\cp = \cm \oplus \ca \oplus \cn$~ with
~$\dim\, \ca=1$.

Let ~$\nu$~ be a (non-unitary) character of ~$A$, ~$\nu\in\ca^*$,
parameterized by a real number ~{\it $d$}, called the {\it conformal
weight} or energy.

Further, let ~ $\mu$ ~ fix a discrete series representation
~$D^\mu$~ of $M$ on the Hilbert space ~$V_\mu\,$, or   the
finite-dimensional (non-unitary) representation of $M$ with the same
Casimirs.

 We call the induced
representation ~$\chi =$ Ind$^G_{P}(\mu\otimes\nu \otimes 1)$~ an
~\blu{\it elementary representation} of $G$ \cite{DMPPT}. (These are
called {\it generalized principal series representations} (or {\it
limits thereof}) in \cite{Knapp}.)   Their spaces of functions are:  $$
\cc_\chi ~=~ \{ \cf \in C^\infty(G,V_\mu) ~ \vr ~ \cf (gman) ~=~
e^{-\nu(H)} \cdot D^\mu(m^{-1})\, \cf (g) \} $$ where ~$a=
\exp(H)\in A'$, ~$H\in\ca'\,$, ~$m\in M'$, ~$n\in N'$. The
representation action is the \blu{left regular action}:  $$
(\ct^\chi(g)\cf) (g') ~=~ \cf (g^{-1}g') ~, \quad g,g'\in G\ .$$
ERs are important due to the following fundamental result:\\
\nt\blu{Theorem \cite{Lan,KnZu}:} Every irreducible admissible representation
of G is equivalent to a subrepresentation of an ER.

\bu An important ingredient in our considerations are the ~\blu{\it
highest/lowest weight representations}~ of ~$\cg^\bac$. These can be
realized as (factor-modules of) Verma modules ~$V^\L$~ over
~$\cg^\bac$, where ~$\L\in (\ch^\bac)^*$, ~$\ch^\bac$ is a Cartan
subalgebra of ~$\cg^\bac$, weight ~$\L = \L(\chi)$~ is determined
uniquely from $\chi$ \cite{Dob}.

Actually, since our ERs may be induced from finite-dimensional
representations of ~$\cm$~ (or their limits) the Verma modules are
always reducible. Thus, it is more convenient to use ~\blu{\it
generalized Verma modules} ~$\tV^\L$~ such that the role of the
highest/lowest weight vector $v_0$ is taken by the
(finite-dimensional) space ~$V_\mu\,v_0\,$. For the generalized
Verma modules (GVMs) the reducibility is controlled only by the
value of the conformal weight $d$. Relatedly, for the \idos{} only
the reducibility w.r.t. non-compact roots is essential.

\bu One main ingredient of our approach is as follows. We group the
(reducible) ERs with the same Casimirs in sets called \red{~{\it
multiplets}} \cite{Dobmul}. The multiplet corresponding to fixed values of the
Casimirs may be depicted as a connected graph, the \blu{vertices} of
which correspond to the reducible ERs and the \blu{lines (arrows)}
between the vertices correspond to intertwining operators.  The
explicit parametrization of the multiplets and of their ERs is
important for understanding of the situation.

In fact, the multiplets contain explicitly all the data necessary to
construct the \idos{}. Actually, the data for each \ido{} consists
of the pair ~$(\b,m)$, where $\b$ is a (non-compact) positive root
of ~$\cg^\bac$, ~$m\in\bbn$, such that the \blu{BGG  Verma module
reducibility condition} (for highest weight modules) is fulfilled \cite{BGG}:
$$ (\L+\r, \b^\vee ) ~=~ m \ , \quad \b^\vee \equiv 2 \b
/(\b,\b) \ $$ $\r$ is half the sum of the positive roots of
~$\cg^\bac$. When the above holds then the Verma module with shifted
weight ~$V^{\L-m\b}$ (or ~$\tV^{\L-m\b}$ ~ for GVM and $\b$
non-compact) is embedded in the Verma module ~$V^{\L}$ (or
~$\tV^{\L}$). This embedding is realized by a singular vector
~$v_s$~ determined by a polynomial ~$\cp_{m,\b}(\cg^-)$~ in the
universal enveloping algebra ~$(U(\cg_-))\ v_0\,$, ~$\cg^-$~ is the
subalgebra of ~$\cg^\bac$ generated by the negative root generators
\cite{Dix}. More explicitly, \cite{Dob}, ~$v^s_{m,\b} = \cp_{m,\b}\, v_0$
(or ~$v^s_{m,\b} = \cp_{m,\b}\, V_\mu\,v_0$ for GVMs). Then
there exists \cite{Dob} an \red{\ido{}} $$  \cd_{m,\b} ~:~ \cc_{\chi(\L)}
~\llr ~ \cc_{\chi(\L-m\b)} $$ given explicitly by: $$
 \cd_{m,\b} ~=~ \cp_{m,\b}(\widehat{\cg^-})  $$ where
~$\widehat{\cg^-}$~ denotes the \blu{right action} on the functions
~$\cf$.

In most of these situations the invariant operator ~$\cd_{m,\b}$~
has a non-trivial invariant kernel in which a subrepresentation of
$\cg$ is realized. Thus, studying the equations with trivial RHS:
$$ \cd_{m,\b}\ f ~=~ 0 \ , \qquad f \in \cc_{\chi(\L)} \ ,$$
is also very important. For example, in many physical applications
 in the case of first order differential operators,
i.e., for ~$m=m_\b = 1$, these equations are called
~\blu{conservation laws}, and the elements ~$f\in \ker \cd_{m,\b}$~
are called ~\blu{conserved currents}.

Below in our exposition we shall use the so-called Dynkin labels: $$
m_i ~\equiv~ (\L+\r,\a^\vee_i)  \ , \quad i=1,\ldots,n, $$ where
~$\L = \L(\chi)$, ~$\r$ is half the sum of the positive roots of
~$\cg^\bac$.

We shall use also   the so-called Harish-Chandra parameters \cite{Har}:
$$ m_\b \equiv (\L+\r, \b )\ ,  $$ where $\b$ is any
positive root of $\cg^\bac$. These parameters are redundant, since
they are expressed in terms of the Dynkin labels, however,   some
statements are best formulated in their terms. (Clearly, both the
Dynkin labels and Harish-Chandra parameters have their origin in the
BGG reducibility condition.)


\section{The  Lie algebra \blu{$su(n,n)$} and parabolically related}

\nt Let ~$\cg ~=~ su(n,n)$, ~$n\geq 2$. The maximal compact subgroup
is ~$\ck \cong u(1)\oplus su(n)\oplus su(n)$, while ~$\cm =
sl(n,\bbc)_\bbr\,$.  The number of ERs in the main
multiplets is equal to \cite{Dobparab} $$\vr W(\cg^\bac,\ch^\bac)\vr\, /\, \vr
W(\cm^\bac,\ch_m^\bac)\vr  =  \left( {2n\atop n}\right)$$ The
signature of the ERs of $\cg$    is: \eqnn{}   \chi  &=&  \{ n_1 ,
\ldots, n_{n-1} , n_{n+1} \ldots, n_{2n-1} ;\ c\, \} \ ,  \quad
n_j \in \bbn\ , \quad c = d - \ha n^2 \eea  The  restricted
Weyl reflection is given by the Knapp--Stein integral operators \cite{KnSt}:
 \eqnn{knast}  && G_{KS}  :  \cc_\chi   \llr   \cc_{\chi'} \
,\qquad
 \chi'  ~=~ \{
(n_1,\ldots,n_{n-1},n_{n+1},\ldots,n_{2n-1})^*  ; \ -c \, \} \ , \\
 &&  (n_1,\ldots,n_{n-1},n_{n+1},\ldots,n_{2n-1})^*  \doteq
 (n_{n+1},\ldots,n_{2n-1},n_1,\ldots,n_{n-1}) \nn\eea

\nt{\bf Multiplets}

 Below we give the multiplets  for  $su(n,n)$ for
 $n=2,3,4$. They are valid also for   ~$sl(2n,\bbr)$~ with ~$\cm$-factor
 ~$sl(n,\bbr) \oplus sl(n,\bbr)$,  and when
 $n=2k$  these are multiplets also for the parabolically related algebra
 ~$su^*(4k)$~ with ~$\cm$-factor $su^*(2k) \oplus su^*(2k)$,

 There are several types of multiplets: the main type, which
contains maximal number of ERs/GVMs, the finite-dimensional and the
discrete series representations, and many reduced types of
multiplets.

The multiplets of the main type are in 1-to-1 correspondence with
the finite-dimensional irreps of  $su(n,n)$, i.e., they will be
labelled by  the  $2n-1$  positive Dynkin labels     $m_i\in\bbn$.

\section{Multiplets of ~SU(2,2), SL(4,$\bbr$) and SU$^*$(4)}

The main multiplet contains six ERs  whose signatures can be given
in the following pair-wise manner: \eqnn{tabl} \chi_0^\pm  &=&  \{ (
m_1, m_3)^\pm  ;  \pm \frac{(m_1 + 2m_2 + m_3)}{2}
   \} \\
\chi'^\pm  &=&  \{  ( m_{12}, m_{23})^\pm  ;  \pm \ha (m_1 + m_3)
   \} \nn\\
\chi''^\pm  &=&  \{  ( m_{2}, m_{13})^\pm  ;  \pm \ha (m_3 - m_1)
    \} \nn \eea
where we have used for the numbers  $m_\b  =  (\L(\chi)+\r,\b)$  the
same compact notation as  for the roots $\b$, and  \eqnn{conutw} &&
(n_1,n_{3})^-  =  (n_1,n_{3})\ , \qquad (n_1,n_{3})^+  =
(n_1,n_{3})^*  =  (n_3,n_{1}) \eea
These multiplets were given first for ~$su^*(4)$~ \cite{DoPe}.

Obviously, the pairs in \eqref{tabl}  are related by Knapp-Stein
integral operators, i.e., \eqn{ackin}  G_{KS}  :  \cc_{\chi^\mp}
\lra \cc_{\chi^\pm} \nn\ee

The multiplets are given explicitly in Fig. 1, where we use the
notation:  $\L^\pm = \L(\chi^\pm)$.  Each \ido\ is represented by an
arrow accompanied by a symbol  $i_{jk}$  encoding the root
 $\a_{jk}$  and the number $m_{\a_{jk}}$ which is involved in the
BGG criterion.  This notation is used to save space, but it can be
used due to the fact that only \idos\ which are non-composite are
displayed, and that the data  $\b,m_\b $, which is involved in the
embedding  $V^\L \lra V^{\L-m_\b,\b}$  turns out to involve only the
 $m_i$  corresponding to simple roots, i.e., for each $\b,m_\b$
there exists  $i = i(\b,m_\b,\L)\in \{ 1,\ldots,2n-1\}$, such that
 $m_\b=m_i $. Hence the data  $\a_{jk} $, $m_{\a_{jk}}$  is
represented by  $i_{jk}$  on the arrows.

\vskip 5mm

 \fig{}{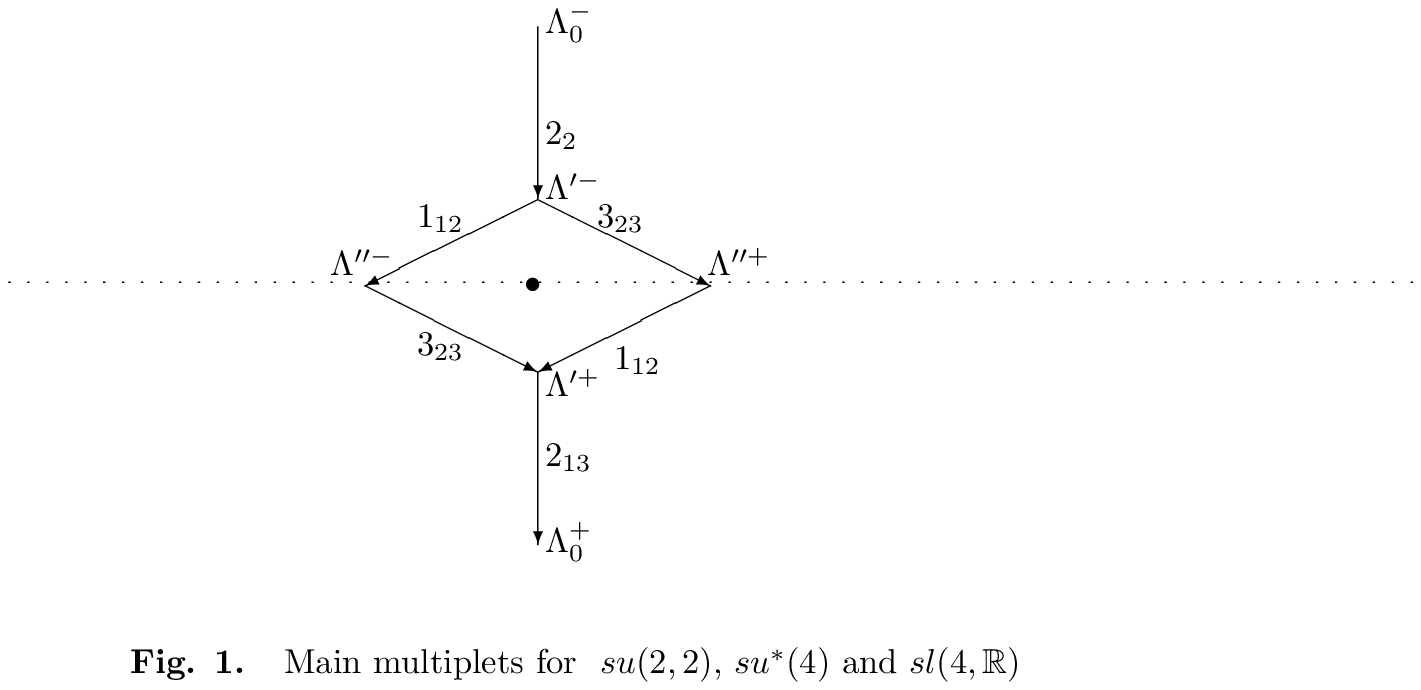}{10cm}

The pairs  $\L^\pm$  are symmetric w.r.t. to the bullet in the
middle of the figure - this represents the Weyl symmetry realized by
the Knapp-Stein operators.

Matters are arranged so that in every multiplet only the ER with
signature  $\chi_0^-$  contains a finite-dimen\-sional nonunitary
subrepresentation in  a finite-dimen\-sional subspace  $\ce$. The
latter corresponds to the finite-dimensional   irrep of  $\cg$  with
signature  $\{ m_1 , m_2 , m_3 \}$  of dimension:  $m_1m_2m_3
m_{12}m_{23}m_{13}/6$. The subspace  $\ce$  is annihilated by the
operator  $G^+ $,\ and is the image of the operator  $G^- $. The
subspace  $\ce$  is annihilated also by the \ido{} acting from
 $\chi^-$  to  $\chi'^-$.
 When all  $m_i=1$  then  $\dim \ce = 1$, and in that case
 $\ce$  is also the trivial one-dimensional UIR of the whole algebra
 $\cg$. Furthermore in that case the conformal weight is zero:
 $d=2+c=2-\ha(m_1+2m_2+m_3)_{\vert_{m_i=1}}=0$.

In the conjugate ER ~$\chi_0^+$~ there is a unitary
subrepresentation in  an infinite-dimen\-sional subspace  $\cd$. It
is annihilated by the operator  $G^- $,\ and is the image of the
operator  $G^+ $.

 All the above is valid also for the algebras ~$sl(4,\bbr) \cong so(3,3)$~ and
 ~$su^*(4) \cong so(5,1)$. However, the latter two do not have
 discrete series representations. On the other hand the algebra
 ~$su(2,2) \cong so(4,2)$~ had discrete series representations and
 furthermore highest/lowest weight series representations.

Thus, in the case of ~$su(2,2)$~  the ER  $\chi_0^+ $ contains both
the  holomorphic discrete series representation and the conjugate
anti-holomorphic discrete series. The direct sum of the latter two
  is realized in the invariant subspace  $\cd$  of the ER  $\chi_0^+ $.

    Note that the corresponding lowest weight GVM
is infinitesimally equivalent only to the holomorphic discrete
series, while the conjugate highest weight GVM is infinitesimally
equivalent to the anti-holomorphic discrete series.\nl The conformal
weight of the ER   $\chi_0^+$  has the restriction  $d = 2+c = 2 +
\ha(m_1+2m_2+m_3)  \geq 4$.

\nt{\bf Remark on SU(1,1)}\nl As we mentioned the case  $su(1,1)$ is
well known - it was studied 60 years ago in the isomorphic form
 $sl(2,\bbr)$  by Gelfand et al \cite{GeNa} and by Bargmann \cite{Barg}.
  In the current setting it was given in \cite{Dobpeds}.
  Here we shall only mention that the multiplets
contain two ERS/GVMs (cf.  ${2n\choose n}_{n=1}=2$), and we can take
as their representatives the pair  $\L^\pm$  and all statements that
fit the setting are true.  In fact, the old results are prototypical
for these pairs, which appear once for each algebra of the conformal
type.\dia

\nt{\bf Reduced multiplets.}\nl There are three types of reduced
multiplets, $R_1 $, $R_2 $, $R_3 $. Each of them contains
 two ERS/GVMs and
may be obtained from the main multiplet by setting formally
 $m_1=0$,  $m_2=0$,  $m_3=0$, resp.  The signatures are \eqnn{redu}
_1\chi^\pm  &=&  \{  ( m_{2}, m_{23})^\pm  ;  \pm \ha m_3
   \}  \\
_2\chi^\pm  &=&  \{  ( m_{1}, m_{3})^\pm  ;  \pm \ha (m_1 + m_3)
   \}  \nn\\
_3  \chi^\pm  &=&  \{  ( m_{12}, m_{2})^\pm  ;  \pm \ha m_1
   \}  \nn \eea

The above is valid for the parabolically related algebras ~$su(2,2), su^*(4), sl(4,\bbr)$.
For ~$su(2,2)$~  the ER  $_2\chi^+$  contains the limits of the
(anti)holomorphic discrete series representations.   Its conformal
weight has the restriction  $d = 2 + \ha (m_1 + m_3) \geq 3$.

Actually, types $R_1 $, $R_3 $  are conjugated under the $^*$
operation (that is not the Weyl symmetry since the sign of  $c$  is
not changed).

Finally, there is the reduced multiplet  $R_{13}$  containing a
single representation \eqn{redus} \chi^s  =  \{  ( m, m)  ;  0  \}
\nn\ee   This multiplet may be omitted from this classification
since it contains no operators, but its importance was understood in
the framework of conformal supersymmetry, i.e., in the multiplet
classification for the superconformal algebra  $su(2,2/N)$ given in
\cite{DPm}. It turns out that the infinite multiplets of
 $su(2,2/N)$  have as building blocks all mentioned above multiplets
of  $su(2,2)$ -   sextets,  doublets  and   singlets.

\section{Multiplets of ~SU(3,3) and SL(6,$\bbr$)}

\nt The main multiplet contains 20 ERs/GVMs whose signatures can be
given in the following pair-wise manner:
\eqnn{tabltri}  \chi_0^\pm
   &=&    \{  (
 m_1,
 m_2,
 m_4,
 m_5)^\pm ; \pm m_\r  \} \\
\chi_a^\pm      &=&      \{  (
 m_1,
 m_{23},
 m_{34},
 m_5)^\pm ; \pm (m_\r  -  m_3)  \} \nn\\
\chi_b^\pm    &=&    \{  (
 m_{12},
 m_{3},
 m_{24},
 m_5)^\pm ; \pm  (m_\r  -  m_{23})  \} \nn\\
\chi_{b'}^\pm    &=&    \{  (
 m_{1},
 m_{24},
 m_{3},
 m_{45})^\pm ; \pm   (m_\r  -  m_{34})  \} \nn\\
\chi_c^\pm      &=&      \{  (
 m_{2},
 m_{3},
 m_{14},
 m_5)^\pm ; \pm (m_\r  -  m_{13})  \} \nn\\
\chi_{c'}^\pm          &=&          \{  (
 m_{12},
 m_{34},
 m_{23},
 m_{45})^\pm ; \pm (m_\r  -  m_{24})  \} \nn\\
\chi_{c''}^\pm      &=&      \{  (
 m_{1},
 m_{25},
 m_{3},
 m_{4})^\pm ; \pm (m_\r  -  m_{35})  \} \nn\\
\chi_d^\pm        &=&        \{  (
 m_{2},
 m_{34},
 m_{13},
 m_{45})^\pm ; \pm (m_\r  -  m_{14})  \} \nn\\
\chi_{d'}^\pm        &=&        \{  (
 m_{12},
 m_{35},
 m_{23},
 m_{4})^\pm ; \pm (m_\r  -  m_{25})  \} \nn\\
\chi_e^\pm      &=&      \{  (
 m_{2},
 m_{35},
 m_{13},
 m_{4})^\pm ; \pm  (m_\r  -  m_{15})    \} \nn  \eea
 where ~$m_\r
 =  \ha(   m_1+ 2m_{2} + 3m_3 + 2m_{4} + m_{5})$.
 They are given in Fig. 2.

\fig{}{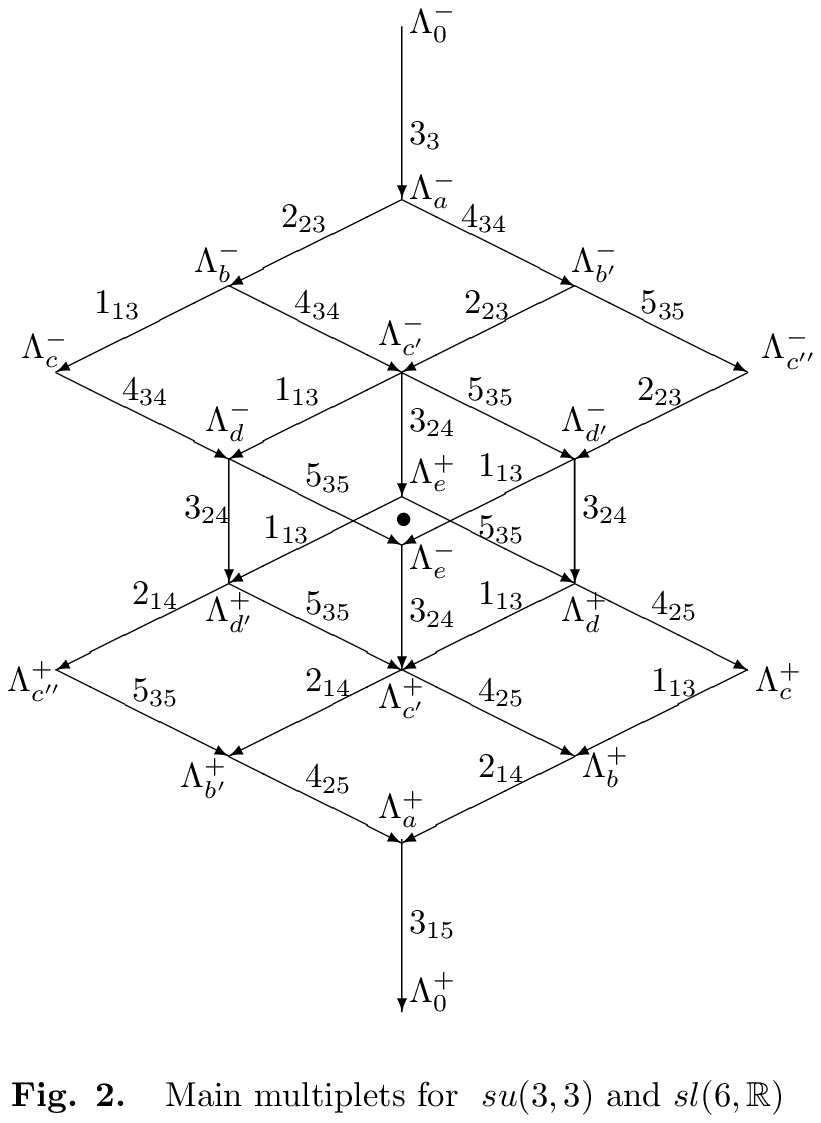}{9cm}

 All general
facts that were stated in the $SU(2,2)$ case are valid also here, in
particular, the special role of the pair
 $\chi^\pm_0 $. The finite-dimensional irreps  $\ce$  of
 $su(3,3)$ or $sl(6,\bbr)$ are sitting in the ERs $\chi^-_0 $ and have dimension
as the UIRs of $SU(6)$.

\bigskip

\nt{\bf Reduced multiplets.}\nl There are five types of reduced
multiplets, $R^3_a $, $a=1,\ldots,5$,  which may be obtained from
the main multiplet by setting formally  $m_a=0$. Multiplets of type
$R^3_4 $, $R^3_5 $, are conjugate  to the multiplets of type $R^3_2
$, $R^3_1 $, resp., and are not shown.

The reduced multiplets of type $R^3_3$ contain 14 ERs/GVMs with
signatures:\eqnn{tabltrtr} \chi_0^\pm ~&=&~ \{\, (
 m_1,
 m_2,
 m_4,
 m_5)^\pm\,;\,\pm m_\r  \,\}  \\
\chi_b^\pm ~&=&~ \{\, (
 m_{12},
 0,
 m_{24},
 m_5)^\pm\,;\,\pm  (m_\r  -  m_{2}) \,\} \cr
\chi_{b'}^\pm ~&=&~ \{\, (
 m_{1},
 m_{24},
 0,
 m_{45})^\pm\,;\,\pm (m_\r  -  m_{4}) \,\} \cr
\chi_c^\pm ~&=&~ \{\, (
 m_{2},
 0,
 m_{14},
 m_5)^\pm\,;\,\pm (m_\r  -  m_{12}) \,\} \cr
\chi_{c''}^\pm ~&=&~ \{\, (
 m_{1},
 m_{25},
 0,
 m_{4})^\pm\,;\,\pm (m_\r  -  m_{45}) \,\} \cr
\chi_d^\pm ~&=&~ \{\, (
 m_{2},
 m_{4},
 m_{12},
 m_{45})^\pm\,;\,\pm (m_\r  -  m_{12,4}) = \mp (m_\r  -  m_{2,45}) \,\} \cr
\chi_e^\pm ~&=&~ \{\, (
 m_{2},
 m_{45},
 m_{12},
 m_{4})^\pm\,;\,\pm (m_\r  -  m_{2,4}) = \pm \ha (m_1+m_5) \,\} \ ,\nn\eea
 here ~$m_\r = \ha (m_1 + 2m_2 + 2m_4 + m_5)$.
  These multiplets are given in Fig. 3.
 They may be called the main type of reduced multiplets since for ~$su(3,3)$~
 in ~$\chi_0^+$~ are contained the limits of the
(anti)holomorphic  discrete series.

\fig{}{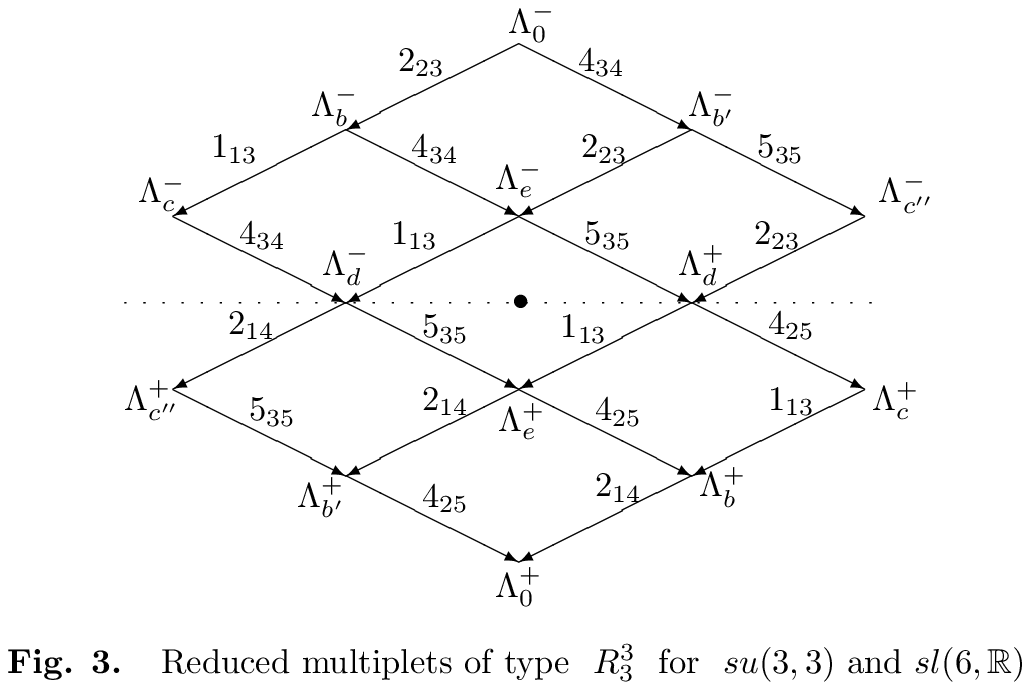}{9cm}

The reduced multiplets of type $R^3_2$, resp., $R^3_1$,  contain 14
ERs/GVMs each. These multiplets are given in Fig. 4, resp., Fig. 5~:

\vskip 5mm

\fig{}{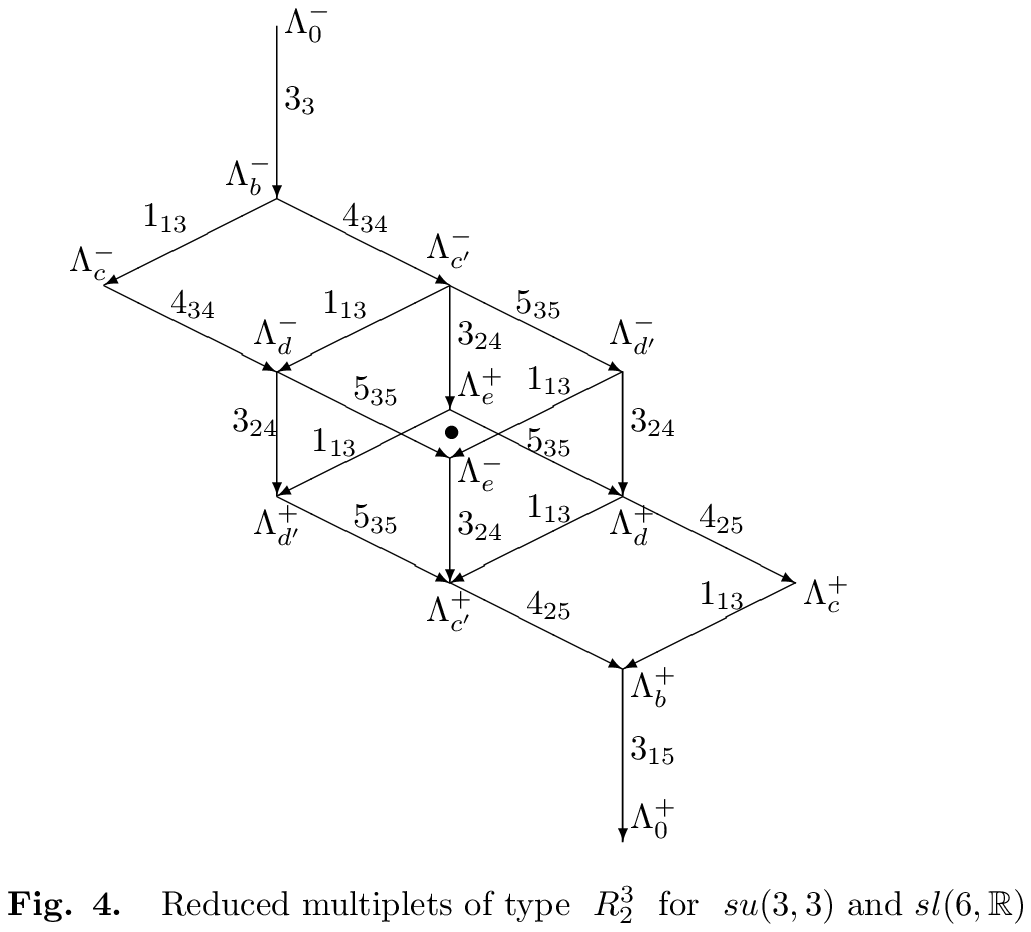}{8cm} \fig{}{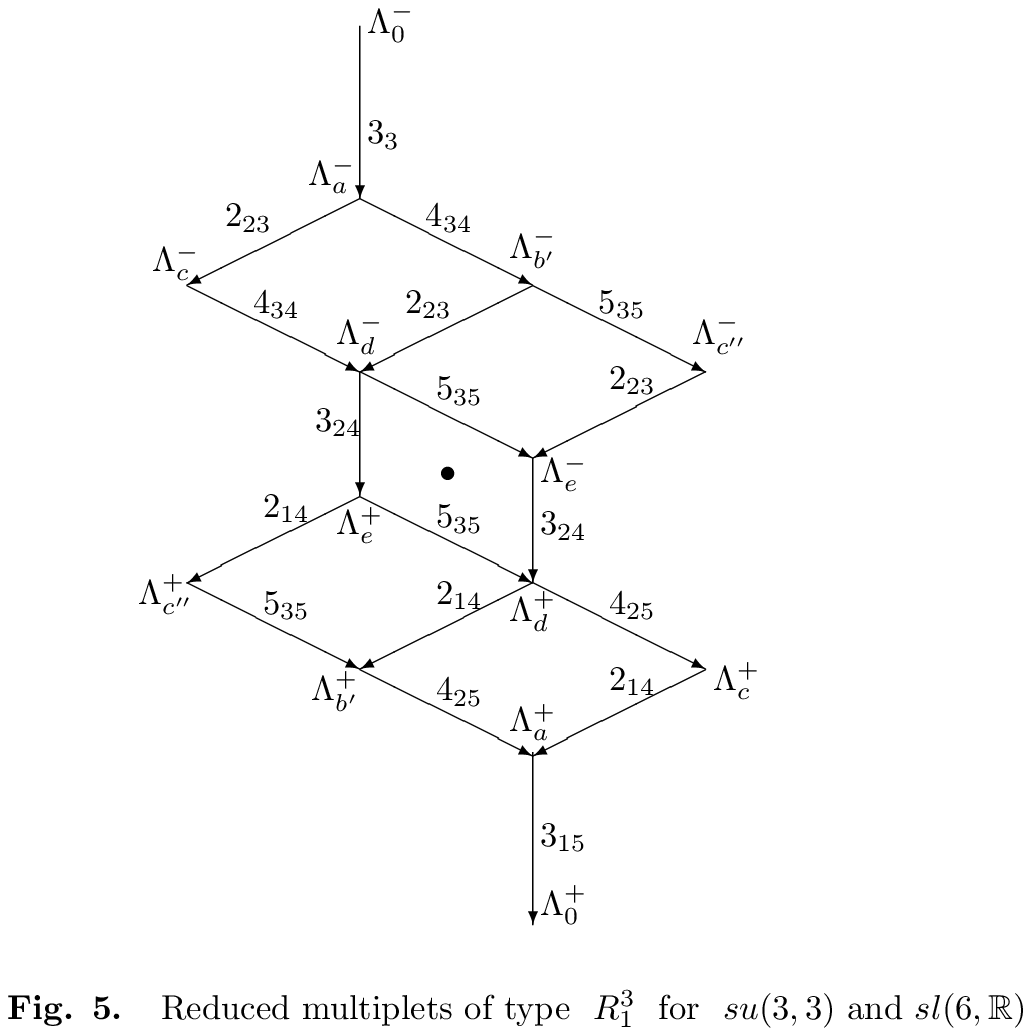}{8cm}

\nt{\bf Further reduction of multiplets}

There are further reductions of the multiplets denoted by
~$R^3_{ab}\,$, $a,b=1,\ldots,5$, $a< b$, which may be obtained from
the main multiplet by setting formally ~$m_a=m_b=0$. From these ten
reductions four (for $(a,b)=(1,2),(2,3),(3,4),(4,5)$) do not
contain representations of physical interest, i.e., induced from
finite-dimensional irreps of the ~$\cm$~ subalgebra. From the others
~$R^3_{35}$~ and ~$R^3_{25}$~ are conjugated to ~$R^3_{13}$~ and
~$R^3_{14}\,$, resp., as explained above. Thus, we present   only
four types of multiplets.

The reduced multiplets of type $R^3_{13}$ contain 10 ERs/GVMs with
signatures:
\eqnn{tabltrona}
\chi_a^\pm ~&=&~ \{\, (
 0,
 m_2,
 m_4,
 m_5)^\pm\,;\,\pm m_\r  \,\} \\
\chi_b^\pm ~&=&~ \{\, (
 m_{2},
0,
 m_{2,4},
 m_5)^\pm\,;\,\pm (m_\r - m_{2}) \,\} \cr
\chi_{b'}^\pm ~&=&~ \{\, (
 0,
 m_{2,4},
 0,
 m_{45})^\pm\,;\,\pm (m_\r - m_{4})\,\} \cr
\chi_{c}^\pm ~&=&~ \{\, (
 0,
 m_{2,45},
 0,
 m_{4})^\pm\,;\,\pm (m_\r - m_{45}) \,\} \cr
\chi_d^\pm ~&=&~ \{\, (
 m_{2},
 m_{4},
 m_{2},
 m_{45})^\pm\,;\,\pm (m_\r - m_{2,4}) = \pm\ha m_5\,\}
 \ ,\nn \eea
  here ~$m_\r = m_2 + m_4 + \ha m_5$.
The multiplets are given in Fig. 6.

\vskip 5mm

\fig{}{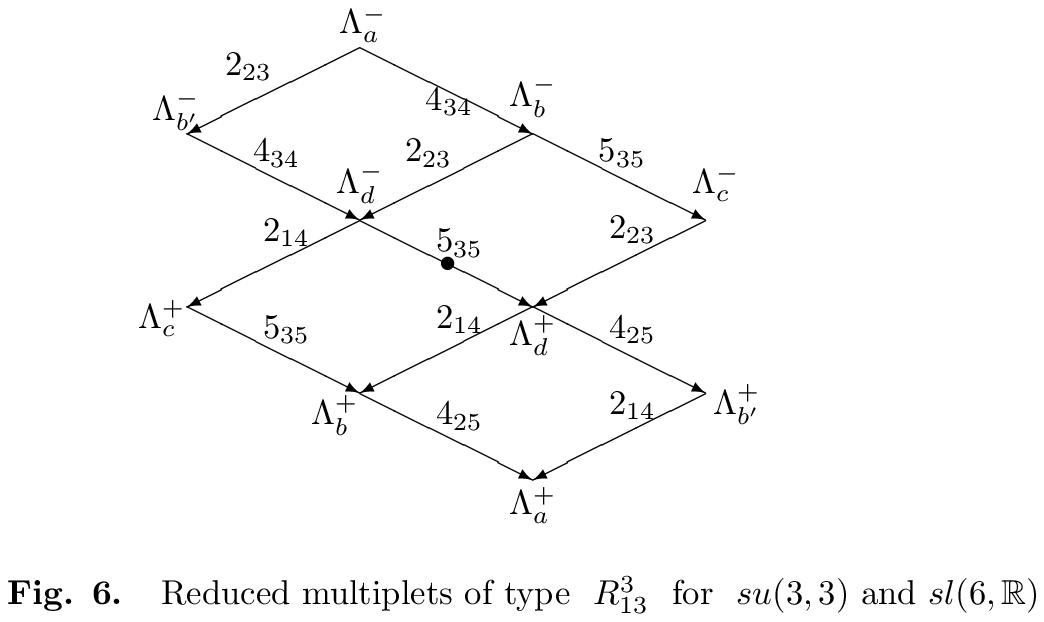}{10cm}

Note that the differential operator (of order $m_5$) from ~$\chi_d^-$~
 to ~$\chi_d^+$~ is a degeneration of an integral Knapp-Stein operator.

The reduced multiplets of type $R^3_{15}$ contain 10 ERs/GVMs with
signatures:
\eqnn{tabltronc}
\chi_0^\pm ~&=&~ \{\, (
 0,
 m_2,
 m_4,
 0)^\pm\,;\,\pm m_\r  \,\} \\
\chi_a^\pm ~&=&~ \{\, (
 0,
 m_{23},
 m_{34},
 0)^\pm\,;\,\pm (m_\r - m_3  )\,\} \cr
\chi_b^\pm ~&=&~ \{\, (
 m_{2},
 m_{3},
 m_{24},
 0)^\pm\,;\,\pm (m_\r - m_{23}) \,\} \cr
\chi_{b'}^\pm ~&=&~ \{\, (
 0,
 m_{24},
 m_{3},
 m_{4})^\pm\,;\,\pm (m_\r - m_{34})\,\} \cr
 \chi_d^\pm ~&=&~ \{\, (
 m_{2},
 m_{34},
 m_{23},
 m_{4})^\pm\,;\,\pm (m_\r - m_{24})=\pm \ha m_3 \,\}
 \ , \nn\eea
  here ~$m_\r = m_2 + \trha m_3 + m_4$.
The multiplets are given in Fig. 7. Here  the differential operator
(of order $m_3$) from ~$\chi_d^-$~
 to ~$\chi_d^+$~  is a degeneration  of an integral Knapp-Stein operator.

\vskip 5mm

\fig{}{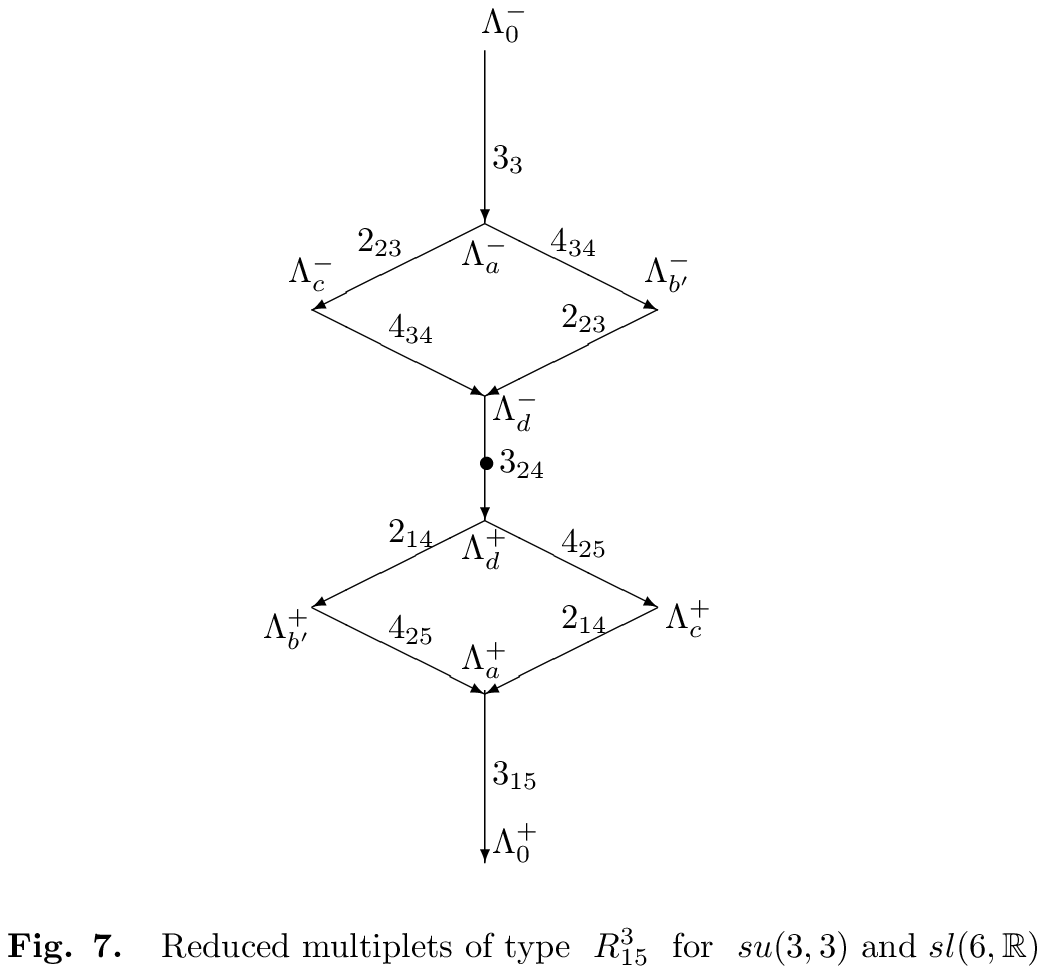}{10cm}

The reduced multiplets of type $R^3_{14}$, $R^3_{24}$,  contain 10 ERs/GVMs each,
the corresponding multiplets being given below:

\fig{}{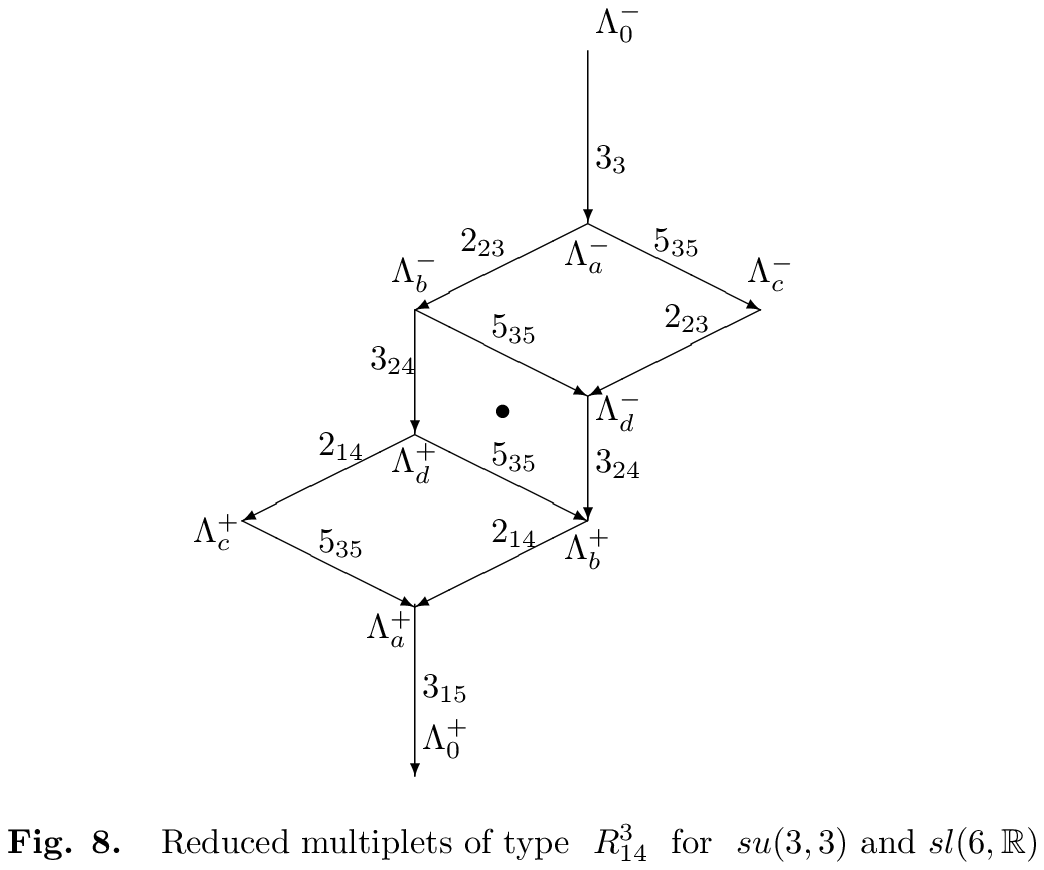}{10cm}
\fig{}{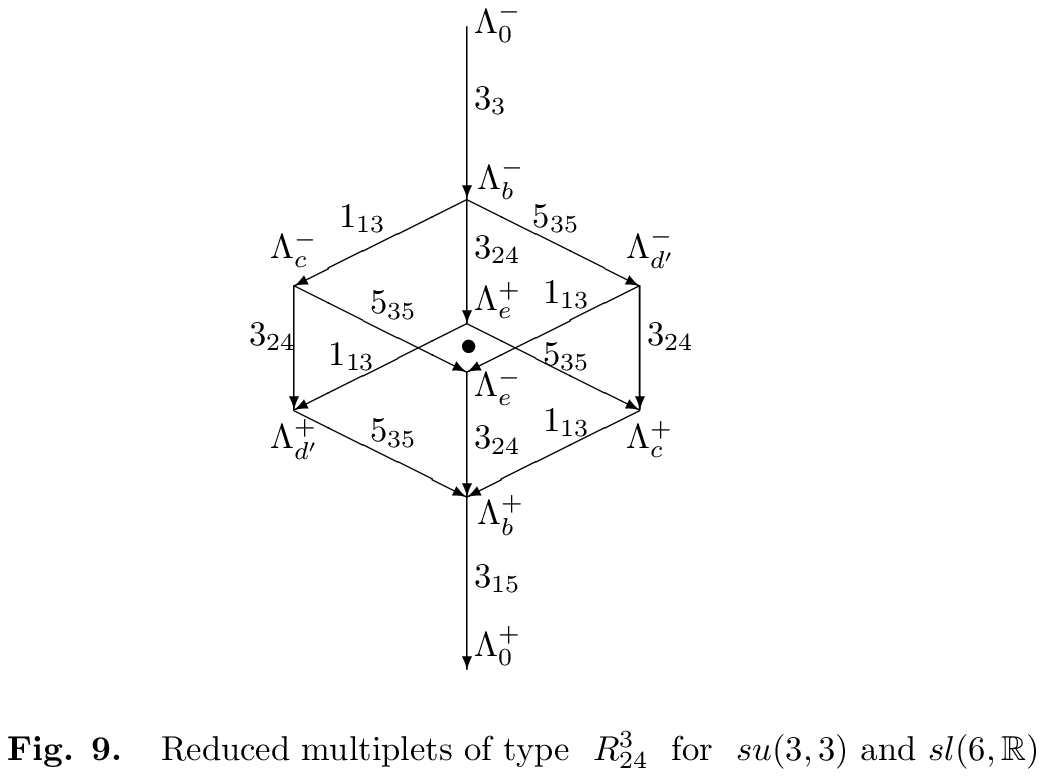}{10cm}

\nt{\bf Last reduction of multiplets}

There are further reductions of the multiplets - triple and quadruple,
but only one triple reduction contains representations of physical interest.
Namely, this is the multiplet  $R^3_{135}\,$,  which
may be obtained from the main multiplet by setting formally  $m_1=m_3=m_5=0$.
It contains 7 ERs/GVMs with signatures:
\eqnn{tabltrona}
&&\chi_a^\pm  =  \{\, (
 0,
 m_2,
 m_4,
0)^\pm\,;\,\pm m_\r = \pm  m_{2,4} \,\} \\
&&\chi_{b}^\pm  =  \{\, (
 0,
 m_{2,4},
 0,
 m_{4})^\pm\,;\,\pm   m_{2}\,\} \nn\\
&&\chi_{b'}^\pm  =  \{\, (
 m_{2},
 0,
 m_{2,4},
0)^\pm\,;\,\pm   m_{4}\,\} \nn\\
&&\chi_d  =  \{\, (
 m_{2},
 m_{4},
 m_{2},
 m_{4})\,;\, 0 \,\} \nn\eea
The multiplets are given below: \fig{}{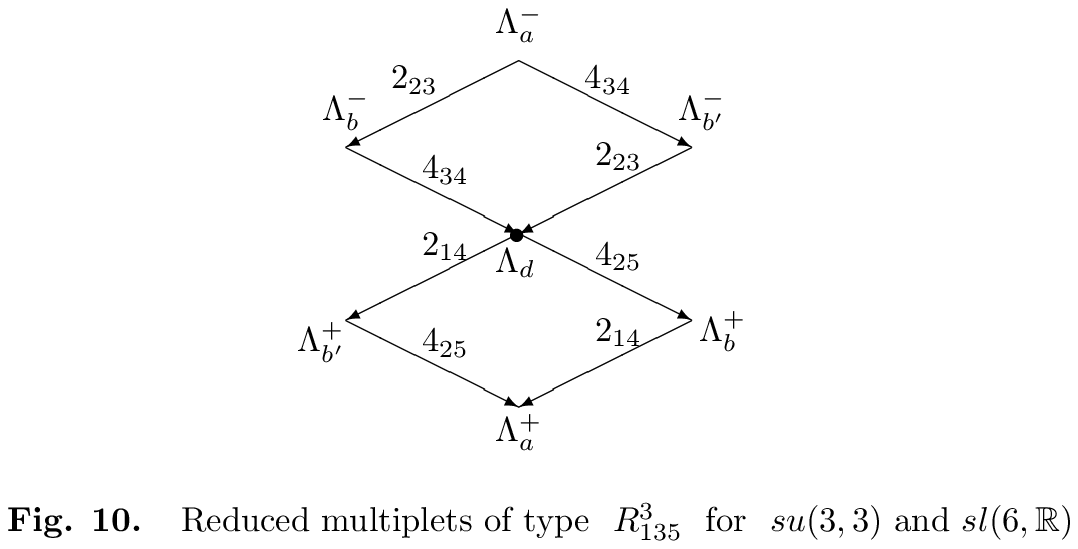}{10cm}
\nt The representation  $\chi^d$  is a singlet, not in a pair, since it has zero weight  $c$,
and the  $\cm$  entries are self-conjugate.  It is placed in the middle
of the figure as the bullet. That ER contains the ~ {\it minimal irreps} ~   characterized by two positive integers
which are denoted in this context as  $m_2\,,m_4\,$. Each such irrep is the kernel of the two
invariant differential operators  $\cd^{m_2}_{14}$  and  $\cd^{m_4}_{25}$, which are of order  $m_2\,$,  $m_4\,$, resp.,
 corresponding to the noncompact roots  $\a_{14}\,$,  $\a_{25}\,$, resp.

\section{Multiplets of ~SU(4,4), SL(8,$\bbr$) and SU$^*$(8)}

The main multiplet ~$R^4$~ contains 70 ERs/GVMs whose
signatures can be given in the following pair-wise manner:
\eqnn{tablf}
&&\chi_0^\pm ~=~ \{\, (
 m_1,
 m_2,
 m_3,
 m_5,
 m_6,
 m_7)^\pm\,;\,\pm  m_\r \,\} \\
&&\chi_{00}^\pm ~=~ \{\, (
 m_1,
 m_2,
 m_{34},
 m_{45},
 m_6,
 m_7)^\pm\,;\,\pm ( m_\r - m_{4}) \,\} \cr
&&\chi_{10}^\pm ~=~ \{\, (
 m_1,
 m_{23},
 m_{4},
 m_{35},
 m_6,
 m_7)^\pm\,;\,\pm ( m_\r - m_{34}) \,\} \cr
 &&\chi_{01}^\pm ~=~ \{\, (
 m_1,
 m_{2},
 m_{35},
 m_{4},
 m_{56},
 m_7)^\pm\,;\,\pm  ( m_\r - m_{45}) \,\} \cr
&&\chi_{20}^\pm ~=~ \{\, (
 m_{12},
 m_{3},
 m_{4},
 m_{25},
 m_6,
 m_7)^\pm\,;\,\pm  (m_\r-m_{24})\,\} \cr
&&\chi_{11}^\pm ~=~ \{\, (
 m_1,
 m_{23},
 m_{45},
 m_{34},
 m_{56},
 m_7)^\pm\,;\,\pm (m_\r -m_{35})\,\} \cr
&&\chi_{02}^\pm ~=~ \{\, (
 m_1,
 m_{2},
 m_{36},
 m_{4},
 m_{5},
 m_{67})^\pm\,;\,\pm  (m_\r-m_{46})\,\} \cr
 &&\chi_{30}^\pm ~=~ \{\, (
 m_{2},
 m_{3},
 m_{4},
 m_{15},
 m_6,
 m_7)^\pm\,;\,\pm  (m_\r  -m_{14} )\,\} \cr
&&\chi_{21}^\pm ~=~ \{\, (
 m_{12},
 m_{3},
 m_{45},
 m_{24},
 m_{56},
 m_7)^\pm\,;\,\pm  (m_\r-m_{25})\,\} \cr
&&\chi_{12}^\pm ~=~ \{\, (
 m_1,
 m_{23},
 m_{46},
 m_{34},
 m_{5},
 m_{67})^\pm\,;\,\pm  (m_\r-m_{36}) \,\} \cr
 &&\chi_{03}^\pm ~=~ \{\, (
 m_1,
 m_{2},
 m_{37},
 m_{4},
 m_{5},
 m_{6})^\pm\,;\,\pm  (m_\r-m_{47})  \,\} \cr
&&\chi_{31}^\pm ~=~ \{\, (
 m_{2},
 m_{3},
 m_{45},
 m_{14},
 m_{56},
 m_7)^\pm\,;\,\pm  (m_\r -m_{15} ) \,\} \cr
&&\chi_{22}^\pm ~=~ \{\, (
 m_{12},
 m_{3},
 m_{46},
 m_{24},
 m_{5},
 m_{67})^\pm\,;\,\pm  (m_\r-m_{26}) \,\} \cr
&&\chi_{13}^\pm ~=~ \{\, (
 m_1,
 m_{23},
 m_{47},
 m_{34},
 m_{5},
 m_{6})^\pm\,;\,\pm   (m_\r-m_{37})\,\} \cr
 &&\chi_{32}^\pm ~=~ \{\, (
 m_{2},
 m_{3},
 m_{46},
 m_{14},
 m_{5},
 m_{67})^\pm\,;\,\pm  (m_\r -m_{16}) \,\} \cr
&&\chi_{23}^\pm ~=~ \{\, (
 m_{12},
 m_{3},
 m_{47},
 m_{24},
 m_{5},
 m_{6})^\pm\,;\,\pm  (m_\r-m_{27}) \,\} \cr
&&\chi_{33}^\pm ~=~ \{\, (
 m_{2},
 m_{3},
 m_{47},
 m_{14},
 m_{5},
 m_{6})^\pm\,;\, \pm  (m_\r-m_{17})    \,\} \cr
&&\chi_{00}'^\pm ~=~ \{\, (
 m_1,
 m_{24},
 m_{5},
 m_{3},
 m_{46},
 m_7)^\pm\,;\,\pm  (m_\r-m_{35}-m_{4}) \,\} \cr
&&\chi_{10}'^\pm ~=~ \{\, (
 m_{12},
 m_{34},
 m_{5},
 m_{23},
 m_{46},
 m_7)^\pm\,;\,\pm  (m_\r-m_{25}-m_{4})  \,\} \cr
&&\chi_{01}'^\pm ~=~ \{\, (
 m_1,
 m_{24},
 m_{56},
 m_{3},
 m_{45},
 m_{67})^\pm\,;\,\pm  (m_\r-m_{36}-m_{4}) \,\} \cr
&&\chi_{20}'^\pm ~=~ \{\, (
 m_{2},
 m_{34},
 m_{5},
 m_{13},
 m_{46},
 m_7)^\pm\,;\,\pm  (m_\r -m_{15}-m_{4})  \,\} \cr
&&\chi_{11}'^\pm ~=~ \{\, (
 m_{12},
 m_{34},
 m_{56},
 m_{23},
 m_{45},
 m_{67})^\pm\,;\,\pm  (m_\r-m_{26}-m_{4})  \,\} \cr
&&\chi_{02}'^\pm ~=~ \{\, (
 m_1,
 m_{24},
 m_{57},
 m_{3},
 m_{45},
 m_{6})^\pm\,;\,\pm  (m_\r-m_{37}-m_{4}) \,\} \cr
&&\chi_{20}''^\pm ~=~ \{\, (
 m_{13},
 m_{4},
 m_{5},
 m_{2},
 m_{36},
 m_7)^\pm\,;\,\pm (m_\r-m_{25}-m_{34})   \,\} \cr
&&\chi_{21}''^\pm ~=~ \{\, (
 m_{13},
 m_{4},
 m_{56},
 m_{2},
 m_{35},
 m_{67})^\pm\,;\,\pm (m_\r-m_{26}-m_{34})  \,\} \cr
&&\chi_{12}''^\pm ~=~ \{\, (
 m_{12},
 m_{35},
 m_{6},
 m_{23},
 m_{4},
 m_{57})^\pm\,;\,\pm (m_\r-m_{26}-m_{45}) \,\} \cr
&&\chi_{02}''^\pm ~=~ \{\, (
 m_{1},
 m_{25},
 m_{6},
 m_{3},
 m_{4},
 m_{57})^\pm\,;\,\pm  (m_\r-m_{36}-m_{45}) \,\} \cr
&&\chi_{30}'^\pm ~=~ \{\, (
 m_{23},
 m_{4},
 m_{5},
 m_{12},
 m_{36},
 m_7)^\pm\,;\,\pm  (m_\r -m_{15}-m_{34}) \,\} \cr
&&\chi_{21}'^\pm ~=~ \{\, (
 m_{2},
 m_{34},
 m_{56},
 m_{13},
 m_{45},
 m_{67})^\pm\,;\,\pm (m_\r -m_{16} - m_{4})  \,\} \cr
&&\chi_{12}'^\pm ~=~ \{\, (
 m_{12},
 m_{34},
 m_{57},
 m_{23},
 m_{45},
 m_{6})^\pm\,;\,\pm  (m_\r-m_{27}-m_{4}) \,\} \cr
&&\chi_{03}'^\pm ~=~ \{\, (
 m_{1},
 m_{25},
 m_{67},
 m_{3},
 m_{4},
 m_{56})^\pm\,;\,\pm (m_\r -m_{37}-m_{45}) \,\} \cr
 &&\chi_{40}'^\pm ~=~ \{\, (
 m_{3},
 m_{4},
 m_{5},
 m_{1},
 m_{26},
 m_7)^\pm\,;\,\pm  (m_\r -m_{15}-m_{24}) \,\} \cr
&&\chi_{31}'^\pm ~=~ \{\, (
 m_{23},
 m_{4},
 m_{56},
 m_{12},
 m_{35},
 m_{67})^\pm\,;\,\pm (m_\r -m_{16}-m_{34})  \,\} \cr
&&\chi_{22}'^\pm ~=~ \{\, (
 m_2,
 m_{34},
 m_{57},
 m_{13},
 m_{45},
 m_{6})^\pm\,;\,\pm (m_\r-m_{17}-m_{4})   \,\} \cr
 &&\chi_{22}''^\pm ~=~ \{\, (
 m_{13},
 m_{4},
 m_{57},
 m_{2},
 m_{35},
 m_{6})^\pm\,;\,\pm (m_\r-m_{27}-m_{34})  \,\} \ , \nn\eea
where ~$m_\r =  \ha(   m_1+ 2m_{2} + 3m_3 + 4m_{4} + 3m_{5}+ 2m_{6}+ m_{7})$.
The multiplets are given explicitly in Fig. 11 (first in
\cite{Dobsunn}).

\fig{}{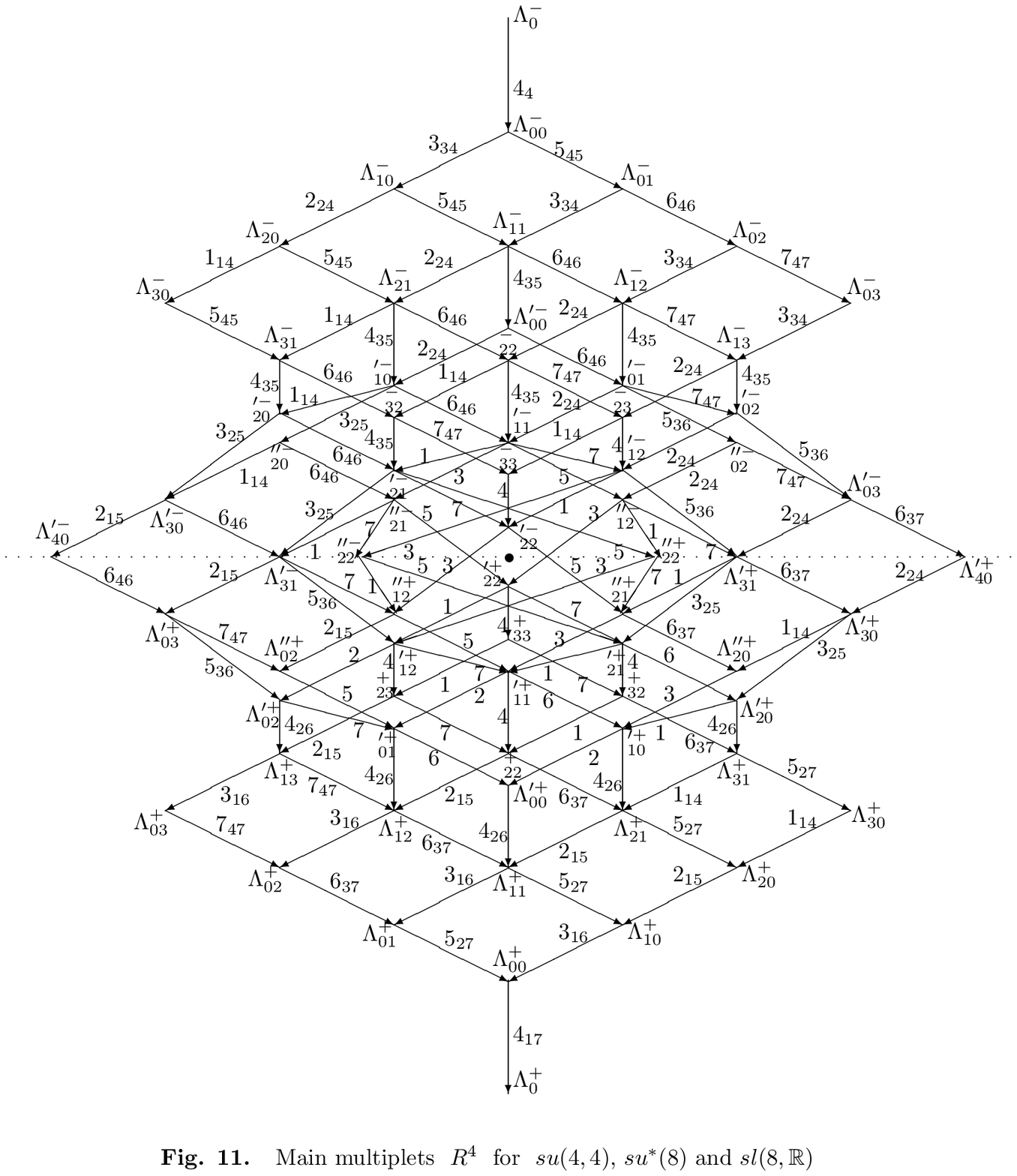}{17cm}

\nt{\bf Main reduced multiplets}

There are nine physically relevant and essentially different reductions of
multiplets denoted by $R^4_3\,$, $R^4_3\,$, $R^4_2\,$, $R^4_1\,$.
 Each of them contains 50 ERs/GVMs   \cite{Dobsuff}.
 Here we present only the  reduced multiplets ~$R^4_4$.
Their 50 ERs/GVMs has
signatures that can be given in the following pair-wise manner:
\eqnn{tablff}
&&\chi_0^\pm ~=~ \{\, (
 m_1,
 m_2,
 m_3,
 m_5,
 m_6,
 m_7)^\pm\,;\,\pm  m_\r \,\}\  \\
 &&\chi_{10}^\pm ~=~ \{\, (
 m_1,
 m_{23},
 0,
 m_{3,5},
 m_6,
 m_7)^\pm\,;\,\pm ( m_\r - m_{3}) \,\} \cr
 &&\chi_{01}^\pm ~=~ \{\, (
 m_1,
 m_{2},
 m_{3,5},
 0,
 m_{56},
 m_7)^\pm\,;\,\pm  ( m_\r - m_{5}) \,\} \cr
&&\chi_{20}^\pm ~=~ \{\, (
 m_{12},
 m_{3},
 0,
 m_{23,5},
 m_6,
 m_7)^\pm\,;\,\pm  (m_\r-m_{23})\,\} \cr
&&\chi_{11}^\pm ~=~ \{\, (
 m_1,
 m_{23},
 m_{5},
 m_{3},
 m_{56},
 m_7)^\pm\,;\,\pm (m_\r -m_{3,5})\,\} \cr
&&\chi_{02}^\pm ~=~ \{\, (
 m_1,
 m_{2},
 m_{3,56},
 0,
 m_{5},
 m_{67})^\pm\,;\,\pm  (m_\r-m_{56})\,\} \cr
 &&\chi_{30}^\pm ~=~ \{\, (
 m_{2},
 m_{3},
 0,
 m_{13,5},
 m_6,
 m_7)^\pm\,;\,\pm  (m_\r  -m_{13} )\,\} \cr
&&\chi_{21}^\pm ~=~ \{\, (
 m_{12},
 m_{3},
 m_{5},
 m_{23},
 m_{56},
 m_7)^\pm\,;\,\pm  (m_\r-m_{23,5})\,\} \cr
&&\chi_{12}^\pm ~=~ \{\, (
 m_1,
 m_{23},
 m_{56},
 m_{3},
 m_{5},
 m_{67})^\pm\,;\,\pm  (m_\r-m_{3,56}) \,\} \cr
 &&\chi_{03}^\pm ~=~ \{\, (
 m_1,
 m_{2},
 m_{3,57},
 0,
 m_{5},
 m_{6})^\pm\,;\,\pm  (m_\r-m_{57})  \,\} \cr
&&\chi_{31}^\pm ~=~ \{\, (
 m_{2},
 m_{3},
 m_{5},
 m_{13},
 m_{56},
 m_7)^\pm\,;\,\pm  (m_\r -m_{13,5} ) \,\} \cr
&&\chi_{22}^\pm ~=~ \{\, (
 m_{12},
 m_{3},
 m_{56},
 m_{23},
 m_{5},
 m_{67})^\pm\,;\,\pm  (m_\r-m_{23,56}) \,\} \cr
&&\chi_{13}^\pm ~=~ \{\, (
 m_1,
 m_{23},
 m_{57},
 m_{3},
 m_{5},
 m_{6})^\pm\,;\,\pm   (m_\r-m_{3,57})\,\} \cr
 &&\chi_{32}^\pm ~=~ \{\, (
 m_{2},
 m_{3},
 m_{56},
 m_{13},
 m_{5},
 m_{67})^\pm\,;\,\pm  (m_\r -m_{13,56}) \,\} \cr
&&\chi_{23}^\pm ~=~ \{\, (
 m_{12},
 m_{3},
 m_{57},
 m_{23},
 m_{5},
 m_{6})^\pm\,;\,\pm  (m_\r-m_{23,57}) \,\} \cr
&&\chi_{33}^\pm ~=~ \{\, (
 m_{2},
 m_{3},
 m_{57},
 m_{13},
 m_{5},
 m_{6})^\pm\,;\, \pm  (m_\r-m_{13,57})    \,\} \cr
&&\chi_{20}''^\pm ~=~ \{\, (
 m_{13},
 0,
 m_{5},
 m_{2},
 m_{3,56},
 m_7)^\pm\,;\,\pm (m_\r-m_{23,5}-m_{3})   \,\} \cr
&&\chi_{21}''^\pm ~=~ \{\, (
 m_{13},
 0,
 m_{56},
 m_{2},
 m_{3,5},
 m_{67})^\pm\,;\,\pm (m_\r-m_{23,56}-m_{3})  \,\} \cr
&&\chi_{12}''^\pm ~=~ \{\, (
 m_{12},
 m_{3,5},
 m_{6},
 m_{23},
 0,
 m_{57})^\pm\,;\,\pm (m_\r-m_{23,56}-m_{5}) \,\} \cr
&&\chi_{02}''^\pm ~=~ \{\, (
 m_{1},
 m_{23,5},
 m_{6},
 m_{3},
 0,
 m_{57})^\pm\,;\,\pm  (m_\r-m_{3,56}-m_{5}) \,\} \cr
&&\chi_{30}'^\pm ~=~ \{\, (
 m_{23},
 0,
 m_{5},
 m_{12},
 m_{3,56},
 m_7)^\pm\,;\,\pm  (m_\r -m_{13,5}-m_{3}) \,\} \cr
&&\chi_{03}'^\pm ~=~ \{\, (
 m_{1},
 m_{23,5},
 m_{67},
 m_{3},
 0,
 m_{56})^\pm\,;\,\pm (m_\r -m_{3,57}-m_{5}) \,\} \cr
 &&\chi_{40}'^\pm ~=~ \{\, (
 m_{3},
 0,
 m_{5},
 m_{1},
 m_{23,56},
 m_7)^\pm\,;\,\pm  (m_\r -m_{13,5}-m_{23}) \,\} \cr
&&\chi_{31}'^\pm ~=~ \{\, (
 m_{23},
 0,
 m_{56},
 m_{12},
 m_{3,5},
 m_{67})^\pm\,;\,\pm (m_\r -m_{13,56}-m_{3})  \,\} \cr
 &&\chi_{22}''^\pm ~=~ \{\, (
 m_{13},
 0,
 m_{57},
 m_{2},
 m_{3,5},
 m_{6})^\pm\,;\,\pm (m_\r-m_{23,57}-m_{3})  \,\}\ , \nn\eea
here ~$m_\r = \ha(   m_1+ 2m_{2} + 3m_3  + 3m_{5}+ 2m_{6}+ m_{7})$.
This is a very important type of reduced
multiplets since for ~$su(4,4)$~ in ~$\chi_0^+$~ are contained the limits of
the (anti)holomorphic  discrete series.
The multiplets are given in Fig. 12.

\fig{}{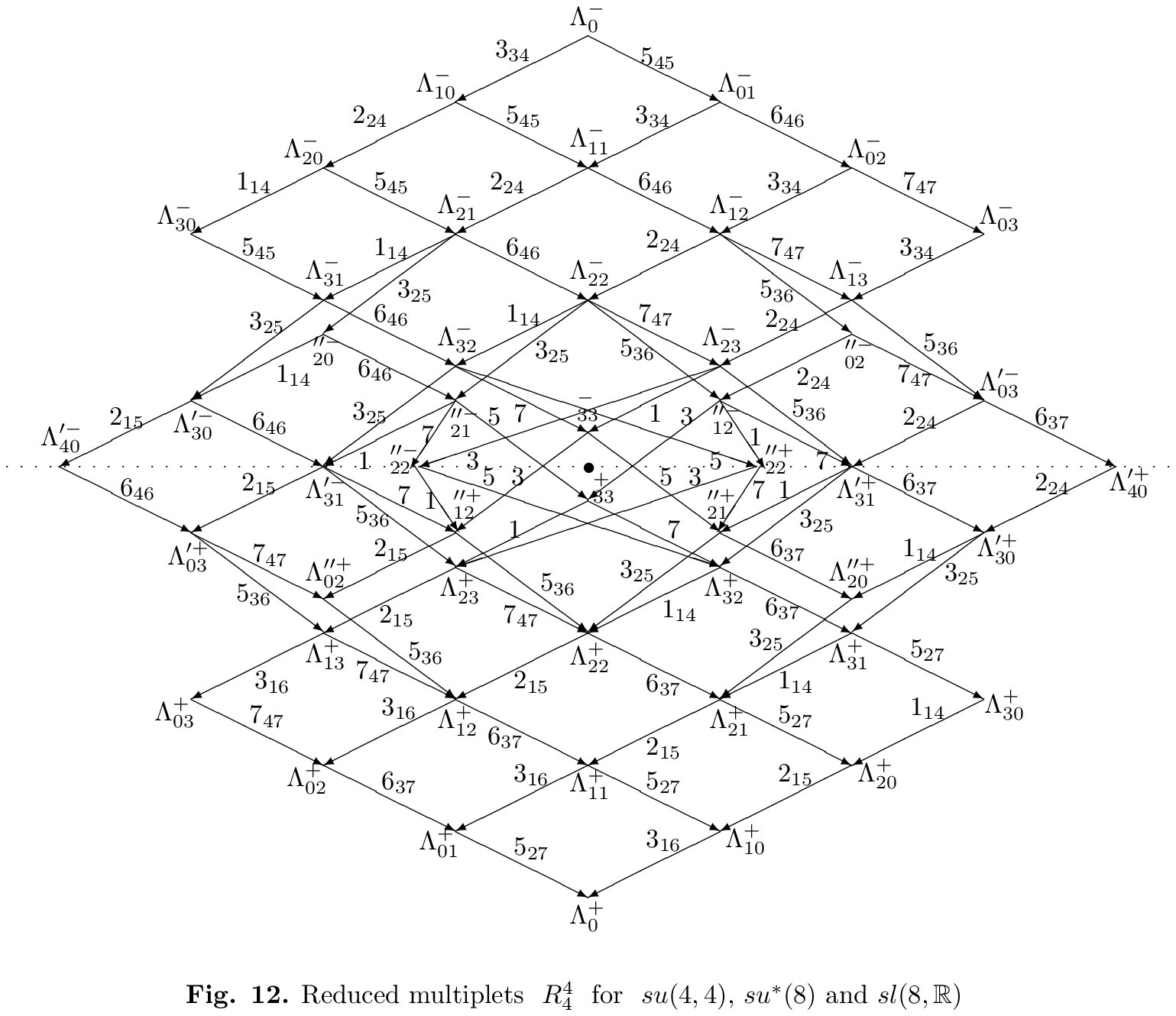}{16cm}

\nt{\bf Further reduction of multiplets}

There are nine physically relevant and essentially different further reductions of
multiplets denoted by ~$R^3_{ab}\,$, $(a,b)=(13),(14),(15),(16),(17),(24),(25),(26),(35)$. They
contain  36 ERs/GVMs each and were given in \cite{Dobsuff}. Here we give only type
$R^4_{13}$~:
\eqnn{tablfot}
&&\chi_0^\pm ~=~ \{\, (
 0,
 m_2,
 0,
 m_5,
 m_6,
 m_7)^\pm\,;\,\pm m_\r  )\,\} \\
&&\chi_{10}^\pm ~=~ \{\, (
 0,
 m_{2},
 m_{4},
 m_{45},
 m_6,
 m_7)^\pm\,;\,\pm ( m_\r - m_{4}) \,\} \cr
&&\chi_{20}^\pm ~=~ \{\, (
 m_{2},
 0,
 m_{4},
 m_{2,45},
 m_6,
 m_7)^\pm\,;\,\pm   (m_\r-m_{2,4}) \,\} \cr
&&\chi_{11}^\pm ~=~ \{\, (
 0,
 m_{2},
 m_{45},
 m_{4},
 m_{56},
 m_7)^\pm\,;\,\pm (m_\r -m_{45}) \,\} \cr
&&\chi_{21}^\pm ~=~ \{\, (
 m_{2},
 0,
 m_{45},
 m_{2,4},
 m_{56},
 m_7)^\pm\,;\,\pm (m_\r-m_{2,45}) \,\} \cr
&&\chi_{12}^\pm ~=~ \{\, (
 0,
 m_{2},
 m_{46},
 m_{4},
 m_{5},
 m_{67})^\pm\,;\,\pm  (m_\r-m_{46}) \,\} \cr
&&\chi_{22}^\pm ~=~ \{\, (
 m_{2},
 0,
 m_{46},
 m_{2,4},
 m_{5},
 m_{67})^\pm\,;\,\pm  (m_\r-m_{2,46}) \,\} \cr
&&\chi_{13}^\pm ~=~ \{\, (
 0,
 m_{2},
 m_{47},
 m_{4},
 m_{5},
 m_{6})^\pm\,;\,\pm    (m_\r-m_{47}) \,\} \cr
&&\chi_{23}^\pm ~=~ \{\, (
 m_{2},
 0,
 m_{47},
 m_{2,4},
 m_{5},
 m_{6})^\pm\,;\,\pm   (m_\r-m_{2,47}) \,\} \cr
&&\chi_{00}'^\pm ~=~ \{\, (
 0,
 m_{2,4},
 m_{5},
 0,
 m_{46},
 m_7)^\pm\,;\,\pm  (m_\r-m_{5}-2m_{4}) \,\} \cr
&&\chi_{10}'^\pm ~=~ \{\, (
 m_{2},
 m_{4},
 m_{5},
 m_{2},
 m_{46},
 m_7)^\pm\,;\,\pm  (m_\r-m_{2,5}-2m_{4}) \,\} \cr
&&\chi_{01}'^\pm ~=~ \{\, (
 0,
 m_{2,4},
 m_{56},
 0,
 m_{45},
 m_{67})^\pm\,;\,\pm (m_\r-m_{56}-2m_{4}) \,\} \cr
&&\chi_{11}'^\pm ~=~ \{\, (
 m_{2},
 m_{4},
 m_{56},
 m_{2},
 m_{45},
 m_{67})^\pm\,;\,\pm  (m_\r-m_{2,56}-2m_{4}) \,\} \cr
&&\chi_{02}'^\pm ~=~ \{\, (
 0,
 m_{2,4},
 m_{57},
 0,
 m_{45},
 m_{6})^\pm\,;\,\pm  (m_\r-m_{57}-2m_{4}) \,\} \cr
&&\chi_{02}''^\pm ~=~ \{\, (
 0,
 m_{2,45},
 m_{6},
 0,
 m_{4},
 m_{57})^\pm\,;\,\pm  (m_\r-m_{6}-2m_{45}) \,\} \cr
&&\chi_{12}'^\pm ~=~ \{\, (
 m_{2},
 m_{4},
 m_{57},
 m_{2},
 m_{45},
 m_{6})^\pm\,;\,\pm  (m_\r-m_{2,57}-2m_{4}) \,\} \cr
&&\chi_{03}'^\pm ~=~ \{\, (
 0,
 m_{2,45},
 m_{67},
 0,
 m_{4},
 m_{56})^\pm\,;\,\pm   (m_\r -m_{67}-2m_{45}) \,\} \cr
 &&\chi_{40}'^\pm ~=~ \{\, (
 0,
 m_{4},
 m_{5},
 0,
 m_{2,46},
 m_7)^\pm\,;\,\pm  (m_\r -m_{5}-2m_{2,4}) \,\}\ ,\nn \eea
  here ~$m_\r = \ha( 2m_{2} +  4m_4  + 3m_{5}+ 2m_{6}+ m_{7})$.
Their diagram is given in Fig. 13.

\vskip 5mm

\fig{}{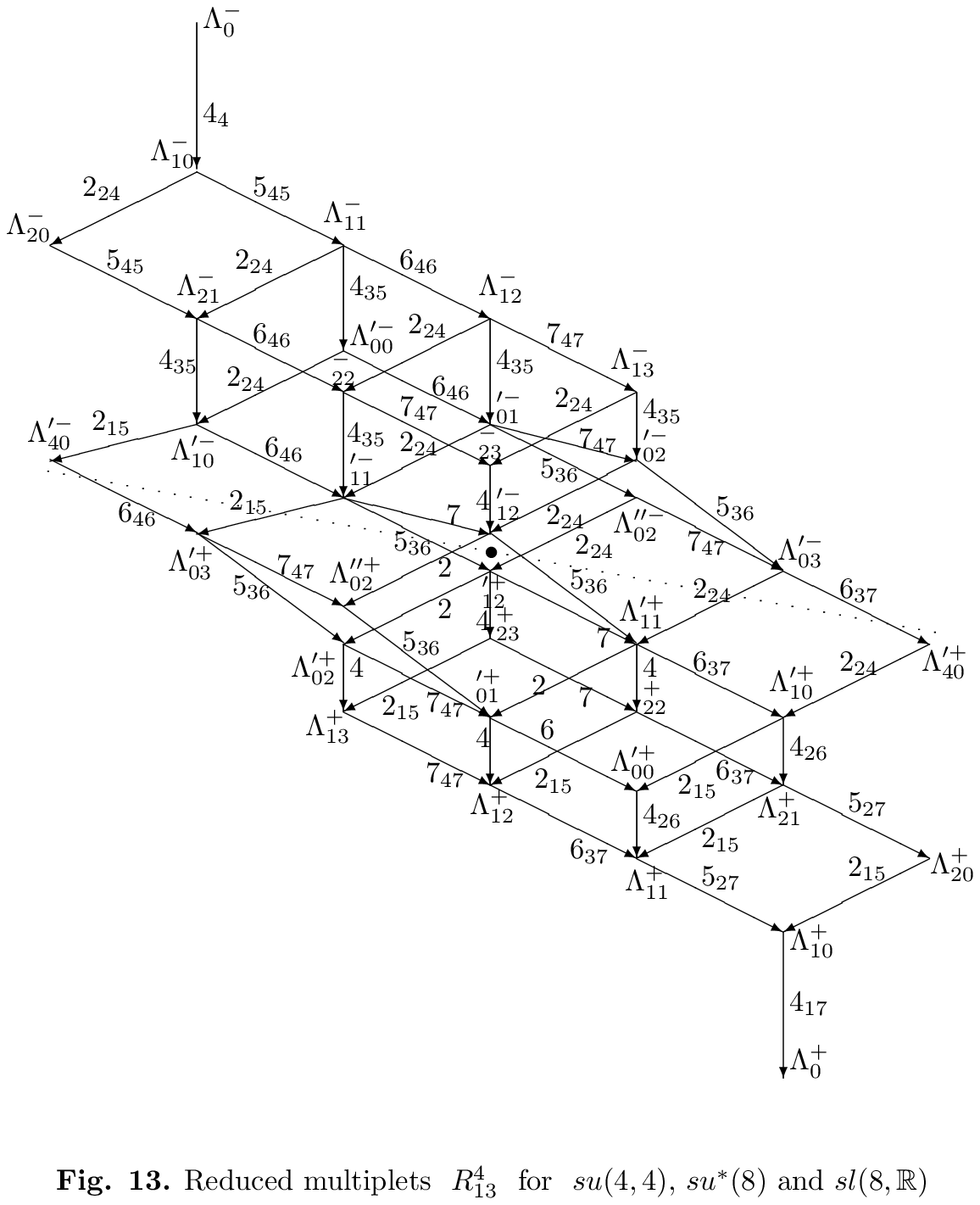}{14cm}

\nt{\bf Yet further reduction of multiplets}

There are six physically relevant and essentially different further reductions of
multiplets denoted by ~$R^3_{abc}\,$, $(a,b,c)=(1,3,5),(1,3,6),(1,3,7),(1,4,6),(1,4,7),(2,4,6)$.
They contain  26 ERs/GVMs each and were given in \cite{Dobsuff}. Here we give only two types.
First we give $R^3_{135}$~:
\eqnn{tablfofit}
&&\chi_0^\pm ~=~ \{\, (
 0,
 m_2,
 0,
 0,
 m_6,
 m_7)^\pm\,;\,\pm m_\r \,\} \\
   &&\chi_{11}^\pm ~=~ \{\, (
 0,
 m_{2},
 m_{4},
 m_{4},
 m_{6},
 m_7)^\pm\,;\,\pm  (m_\r -m_{4}) \,\} \cr
 &&\chi_{21}^\pm ~=~ \{\, (
 m_{2},
 0,
 m_{4},
 m_{2,4},
 m_{6},
 m_7)^\pm\,;\,\pm (m_\r-m_{2,4})  \,\} \cr
&&\chi_{12}^\pm ~=~ \{\, (
 0,
 m_{2},
 m_{4,6},
 m_{4},
 0,
 m_{67})^\pm\,;\,\pm  (m_\r-m_{4,6}) \,\} \cr
  &&\chi_{22}^\pm ~=~ \{\, (
 m_{2},
 0,
 m_{4,6},
 m_{2,4},
 0,
 m_{67})^\pm\,;\,\pm  (m_\r-m_{2,4,6}) \,\} \cr
&&\chi_{13}^\pm ~=~ \{\, (
 0,
 m_{2},
 m_{4,67},
 m_{4},
 0,
 m_{6})^\pm\,;\,\pm   (m_\r-m_{4,67}) \,\} \cr
&&\chi_{23}^\pm ~=~ \{\, (
 m_{2},
 0,
 m_{4,67},
 m_{2,4},
 0,
 m_{6})^\pm\,;\,\pm   (m_\r-m_{2,4,67}) \,\} \cr
&&\chi_{00}'^\pm ~=~ \{\, (
 0,
 m_{2,4},
 0,
 0,
 m_{4,6},
 m_7)^\pm\,;\,\pm  (m_\r -2m_{4}) \,\} \cr
 &&\chi_{01}'^\pm ~=~ \{\, (
 0,
 m_{2,4},
 m_{6},
 0,
 m_{4},
 m_{67})^\pm\,;\,\pm  (m_\r-m_{6}-2m_{4}) \,\} \cr
 &&\chi_{02}'^\pm ~=~ \{\, (
 0,
 m_{2,4},
 m_{67},
 0,
 m_{4},
 m_{6})^\pm\,;\,\pm (m_\r-m_{67}-2m_{4}) \,\} \cr
&&\chi_{30}'^\pm ~=~ \{\, (
 m_{2},
 m_{4},
 0,
 m_{2},
 m_{4,6},
 m_7)^\pm\,;\,\pm (m_\r -m_{2}-2m_{4}) \,\} \cr
&&\chi_{40}'^\pm ~=~ \{\, (
 0,
 m_{4},
 0,
 0,
 m_{26},
 m_7)^\pm\,;\,\pm  (m_\r -2m_{2,4}) \,\} \cr
&&\chi_{31}'^\pm ~=~ \{\, (
 m_{2},
 m_{4},
 m_{6},
 m_{2},
 m_{4},
 m_{67})^\pm\,;\,\pm  (m_\r -m_{2,6}-2m_{4}) = \pm\ha m_7 \,\} \ ,\nn \eea
  here ~$m_\r =  m_{2} +  2m_4 + m_{6}+ \ha m_{7}\,$.
   The multiplets are given in Fig. 3-135.
Note that the differential operator (of order $m_7$) from ~$\chi_{31}^-$~
 to ~$\chi_{31}^+$~ is a degeneration of an integral Knapp-Stein operator.

\vskip 5mm

\fig{}{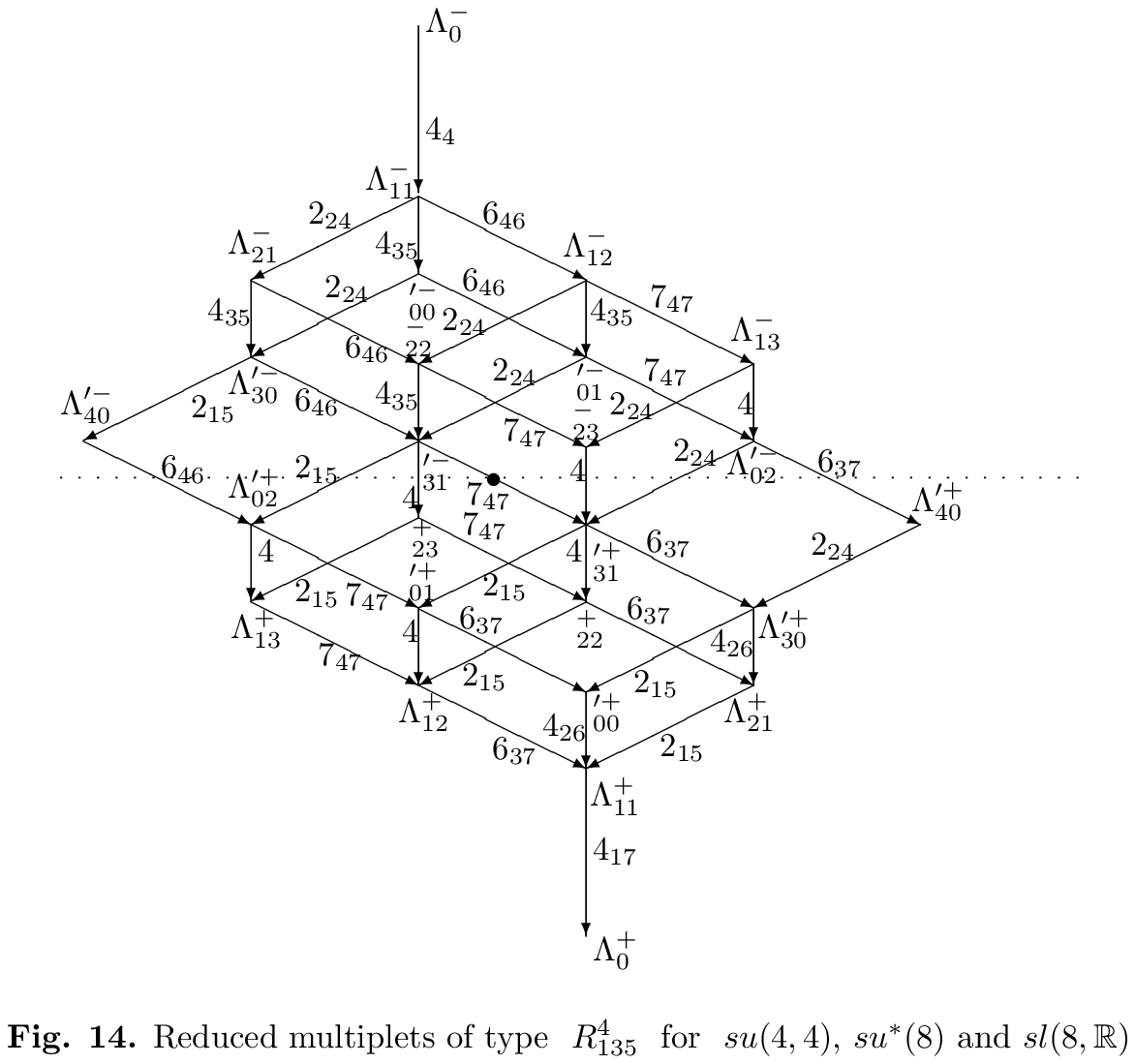}{14cm}

Then we give type ~$R^4_{137}$~:
\eqnn{tablfosetr}
&&\chi_0^\pm ~=~ \{\, (
 0,
 m_2,
 0,
 m_5,
 m_6,
 0)^\pm\,;\,\pm m_\r \,\} \\
 &&\chi_{10}^\pm ~=~ \{\, (
 0,
 m_{2},
 m_{4},
 m_{45},
 m_6,
 0)^\pm\,;\,\pm  ( m_\r - m_{4}) \,\} \cr
  &&\chi_{20}^\pm ~=~ \{\, (
 m_{2},
 0,
 m_{4},
 m_{2,45},
 m_6,
 0)^\pm\,;\,\pm  (m_\r-m_{2,4}) \,\} \cr
&&\chi_{11}^\pm ~=~ \{\, (
 0,
 m_{2},
 m_{45},
 m_{4},
 m_{56},
 0)^\pm\,;\,\pm  (m_\r -m_{45}) \,\} \cr
 &&\chi_{21}^\pm ~=~ \{\, (
 m_{2},
 0,
 m_{45},
 m_{2,4},
 m_{56},
 0)^\pm\,;\,\pm  (m_\r-m_{2,45}) \,\} \cr
&&\chi_{12}^\pm ~=~ \{\, (
 0,
 m_{2},
 m_{46},
 m_{4},
 m_{5},
 m_{6})^\pm\,;\,\pm  (m_\r-m_{46}) \,\} \cr
  &&\chi_{22}^\pm ~=~ \{\, (
 m_{2},
 0,
 m_{46},
 m_{2,4},
 m_{5},
 m_{6})^\pm\,;\, (m_\r-m_{2,46}) \,\} \cr
 &&\chi_{00}'^\pm ~=~ \{\, (
 0,
 m_{2,4},
 m_{5},
 0,
 m_{46},
 0)^\pm\,;\,\pm  (m_\r-m_{5}-2m_{4}) \,\} \cr
 &&\chi_{01}'^\pm ~=~ \{\, (
 0,
 m_{2,4},
 m_{56},
 0,
 m_{45},
 m_{6})^\pm\,;\,\pm  (m_\r-m_{56}-2m_{4}) \,\} \cr
  &&\chi_{30}'^\pm ~=~ \{\, (
 m_{2},
 m_{4},
 m_{5},
 m_{2},
 m_{46},
 0)^\pm\,;\,\pm (m_\r -m_{2,5}-2m_{4}) \,\} \cr
 &&\chi_{03}'^\pm ~=~ \{\, (
 0,
 m_{2,45},
 m_{6},
 0,
 m_{4},
 m_{56})^\pm\,;\,\pm   (m_\r -m_{6}-2m_{45}) \,\} \cr
 &&\chi_{40}'^\pm ~=~ \{\, (
 0,
 m_{4},
 m_{5},
 0,
 m_{2,46},
 0)^\pm\,;\,\pm  (m_\r -m_{5}-2m_{2,4}) \,\} \cr
&&\chi_{31}'^\pm ~=~ \{\, (
 m_{2},
 m_{4},
 m_{56},
 m_{2},
 m_{45},
 m_{6})^\pm\,;\, (m_\r -m_{2,56}-2m_{4})= \pm\ha m_5 \,\} \ ,\nn \eea
  here ~$m_\r = m_{2} +  2m_4  + \trha m_{5}+ m_{6}$.
The multiplets are given in Fig. 15.
Note that the differential operator (of order $m_5$) from ~$\chi_{31}^-$~
 to ~$\chi_{31}^+$~ is a degeneration of an integral Knapp-Stein operator.

\vskip 5mm

\fig{}{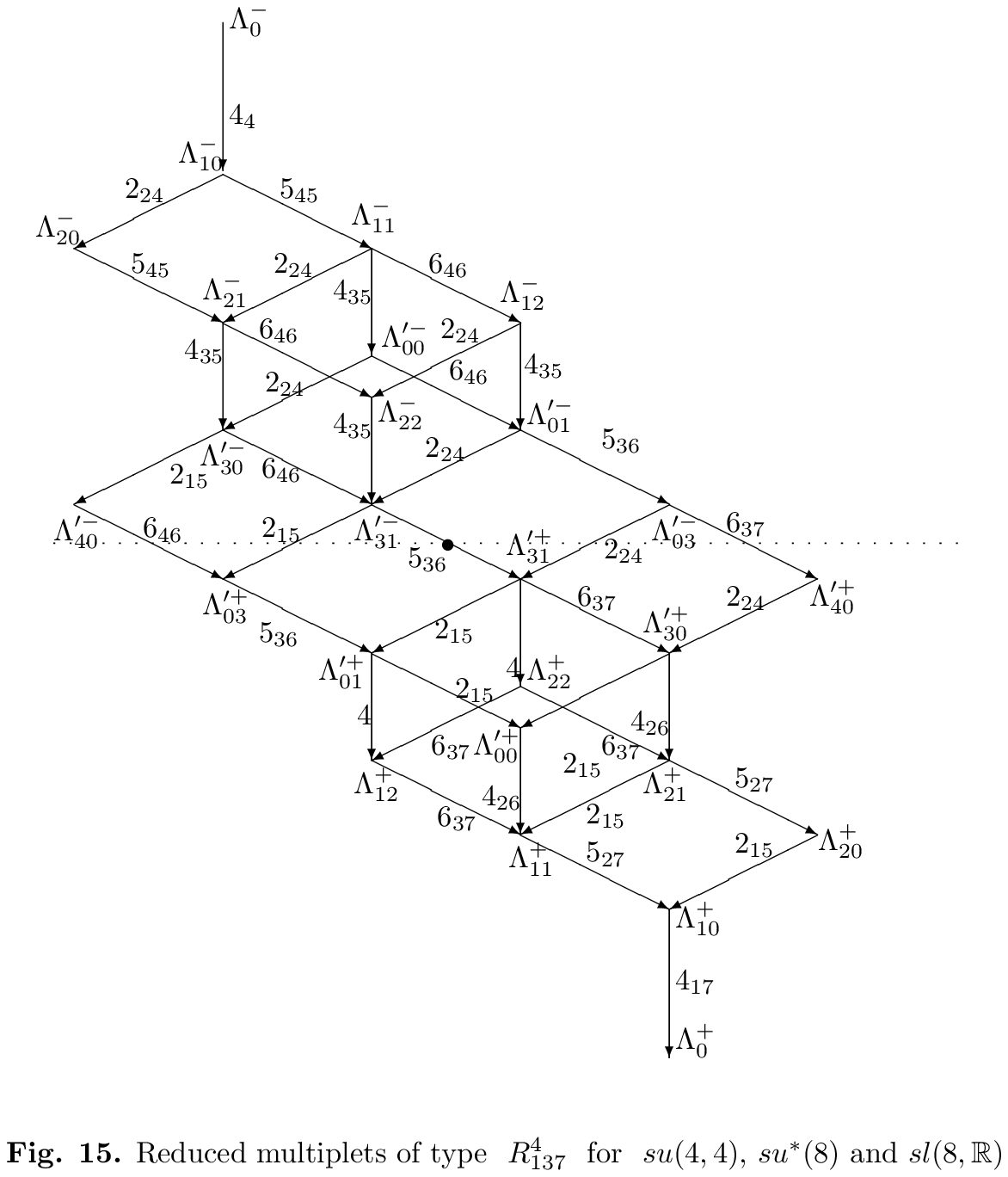}{12cm}

\nt{\bf Last reduction of multiplets}

There are further reductions of the multiplets - quadruple, etc.,
but only one quadruple reduction contains representations of physical interest.
Namely, this is the multiplet ~$R^4_{1357}\,$,  which
may be obtained from the main multiplet by setting formally ~$m_1=m_3=m_5=m_7=0$.
These   multiplets contain 19 ERs/GVMs
whose signatures can be given in the following  manner:
\eqnn{tablfosetrfi}
&&\chi_0^\pm ~=~ \{\, (
 0,
 m_2,
 0,
 0,
 m_6,
 0)^\pm\,;\,\pm m_\r = \pm (m_{2} +  2m_4 + m_{6})\,\} \\
     &&\chi_{11}^\pm ~=~ \{\, (
 0,
 m_{2},
 m_{4},
 m_{4},
 m_{6},
 0)^\pm\,;\,\pm  (m_\r -m_{4}) = \pm m_{2,4,6} \,\} \cr
 &&\chi_{21}^\pm ~=~ \{\, (
 m_{2},
 0,
 m_{4},
 m_{2,4},
 m_{6},
 0)^\pm\,;\,\pm   (m_\r-m_{2,4})= \pm m_{4,6} \,\} \cr
&&\chi_{12}^\pm ~=~ \{\, (
 0,
 m_{2},
 m_{46},
 m_{4},
 0,
 m_{6})^\pm\,;\,\pm  (m_\r-m_{4,6})= \pm m_{2,4} \,\} \cr
  &&\chi_{22}^\pm ~=~ \{\, (
 m_{2},
 0,
 m_{46},
 m_{2,4},
 0,
 m_{6})^\pm\,;\, \pm (m_\r-m_{2,4,6}) = \pm m_4 \,\} \cr
 &&\chi_{00}'^\pm ~=~ \{\, (
 0,
 m_{2,4},
 0,
 0,
 m_{46},
 0)^\pm\,;\,\pm  (m_\r -2m_{4})= \pm m_{2,6}\,\} \cr
 &&\chi_{01}'^\pm ~=~ \{\, (
 0,
 m_{2,4},
 m_{6},
 0,
 m_{4},
 m_{6})^\pm\,;\, \pm (m_\r-m_{6}-2m_{4}) = \pm  m_{2}\,\} \cr
  &&\chi_{30}'^\pm ~=~ \{\, (
 m_{2},
 m_{4},
 0,
 m_{2},
 m_{46},
 0)^\pm\,;\, \pm (m_\r -m_{2}-2m_{4})=\pm  m_{6}\,\} \cr
   &&\chi_{40}'^\pm ~=~ \{\, (
 0,
 m_{4},
 0,
 0,
 m_{2,46},
 0)^\pm\,;\,\pm (m_\r -2m_{2,4})= \pm (m_{6} -m_{2} )\,\} \cr
&&\chi_{31} ~=~ \{\, (
 m_{2},
 m_{4},
 m_{6},
 m_{2},
 m_{4},
 m_{6})\,;\,  0 \,\} \nn \eea
 The multiplets are given in Fig. 3-1357:

\vskip 5mm

 \fig{}{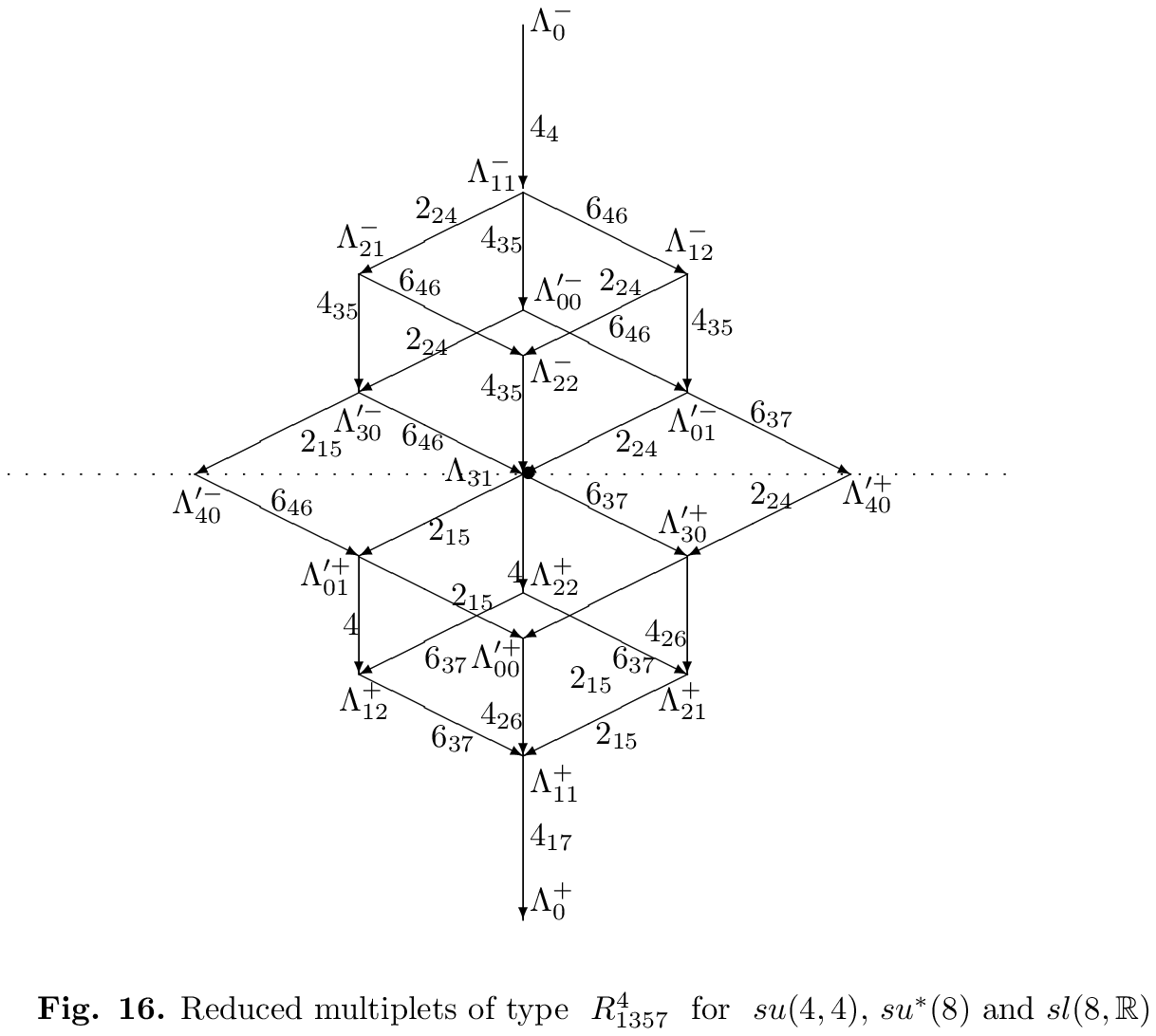}{10cm}

Note that the ER  $\chi_{31}$  is not in a pair -- it has $c=0$ and
its ~$\cm$~ signature is self conjugated. It is placed in the middle
of the figure as the bullet. That ER contains the  ~{\it minimal
irreps}~   characterized by three positive integers
which are denoted in this context as  $m_2 ,m_4 ,m_6 $. Each such
irrep is the kernel of the three invariant differential operators
$\cd^{m_2}_{15} $,  $\cd^{m_4}_{26} $,  $\cd^{m_6}_{37}$, which are
of order  $m_2 $,  $m_4 $,  $m_6 $, resp., and correspond to the
noncompact roots $\a_{15} $,  $\a_{26} $,  $\a_{37} $, resp.

\section*{Acknowledgments} The author thanks the Organizers for the kind
invitation to give a plenary talk  at the International Conference
on Integrable Systems and Quantum Symmetries, Prague,  June
2014.

\end{document}